\documentclass[noname,nonblindrev]{informs-noname}
\OneAndAHalfSpacedXII

\usepackage{natbib}
 \bibpunct[, ]{(}{)}{,}{a}{}{,}%
 \def\bibfont{\small}%

\usepackage{algorithm}
\usepackage{algpseudocode}
\usepackage{makebox}
\usepackage{mathtools}
\usepackage{pgfplotstable}
\usepackage{pgfplots}
\usepackage{pifont}
\usepackage{relsize}
\usepackage[aboveskip=-0.25ex]{subcaption}
\usepackage{array,multirow,graphicx}
\usepackage{tikz}

\usepgfplotslibrary{colorbrewer}
\pgfplotsset{cycle list/Set3-12,compat=1.16}
\pgfplotsset{compat=1.16,cycle list/Dark2-8}
\usetikzlibrary{pgfplots.colorbrewer,pgfplots.fillbetween}

\algnewcommand\algorithmicinput{\textbf{Input:}}
\algnewcommand\Input{\item[\algorithmicinput]}
\algnewcommand\algorithmicoutput{\textbf{Output:}}
\algnewcommand\Output{\item[\algorithmicoutput]}

\pgfplotsset{
    layers/legend behind plots/.define layer set={
            axis background,axis grid,axis ticks,axis lines,axis tick labels,main,axis descriptions,axis foreground
    }{
        grid style= {/pgfplots/on layer=axis grid},
        tick style= {/pgfplots/on layer=axis ticks},
        axis line style= {/pgfplots/on layer=axis lines},
        label style= {/pgfplots/on layer=axis descriptions},
        legend style= {/pgfplots/on layer=axis tick labels}, 
        title style= {/pgfplots/on layer=axis descriptions},
        colorbar style= {/pgfplots/on layer=axis descriptions},
        ticklabel style= {/pgfplots/on layer=axis tick labels},
        axis background@ style={/pgfplots/on layer=axis background},
        3d box foreground style={/pgfplots/on layer=axis foreground},
    },
}

\TheoremsNumberedThrough
\EquationsNumberedThrough

\begin{document}
\RUNAUTHOR{Tasseff et al.}
\TITLE{Polyhedral Relaxations for Optimal Pump Scheduling of Potable Water Distribution Networks}
\RUNTITLE{Polyhedral Relaxations for Optimal Pump Scheduling of Water Networks}
\MANUSCRIPTNO{00-0000-0000.00}

\ARTICLEAUTHORS{%
\AUTHOR{Byron Tasseff, Russell Bent, Carleton Coffrin}
\AFF{Los Alamos National Laboratory}
\AUTHOR{Clayton Barrows, Devon Sigler, Jonathan Stickel, Ahmed S. Zamzam}
\AFF{National Renewable Energy Laboratory}
\AUTHOR{Yang Liu}
\AFF{Stanford University}
\AUTHOR{Pascal Van Hentenryck}
\AFF{Georgia Institute of Technology}
}

\ABSTRACT{%
The classic pump scheduling or Optimal Water Flow (OWF) problem for water distribution networks (WDNs) minimizes the cost of power consumption for a given WDN over a fixed time horizon.
In its exact form, the OWF is a computationally challenging mixed-integer nonlinear program (MINLP).
It is complicated by nonlinear equality constraints that model network physics, discrete variables that model operational controls, and intertemporal constraints that model changes to storage devices.
To address the computational challenges of the OWF, this paper develops tight polyhedral relaxations of the original MINLP, derives novel valid inequalities (or cuts) using duality theory, and implements novel optimization-based bound tightening and cut generation procedures.
The efficacy of each new method is rigorously evaluated by measuring empirical improvements in OWF primal and dual bounds over forty-five literature instances.
The evaluation suggests that our relaxation improvements, model strengthening techniques, and a thoughtfully selected polyhedral relaxation partitioning scheme can substantially improve OWF primal and dual bounds, especially when compared with similar relaxation-based techniques that do not leverage these new methods.%
}


\maketitle

\section{Introduction}
\label{section:introduction}
There are around $52{,}000$ community water systems in the United States, defined as serving twenty-five or more year-round residents.
Nearly $85\%$ of the U.S. population is served by only $5\%$ of these systems, while the remaining $95\%$ comprise small systems, each of which is defined as serving $3{,}300$ people or fewer.
Since they often lack economies of scale, these smaller utilities pay more per unit of water produced than larger utilities.
Furthermore, in these small systems, nearly all of the energy consumed is via electricity, and $80\%$ is consumed by motors used for pumping operations \citep{copeland2014energy}.
Pumping is thus a primary operational cost.
Additionally, since these small systems typically lack smart controls and computational infrastructure, they often rely on ad hoc rules at the expense of system efficiency \citep{salomons2020practical}.
As such, the ability to more cost-effectively control pumps for small water distribution networks (WDNs) could benefit the reliability and overall efficiency of water infrastructure in the United States and beyond.

This goal is well-posed by the classic pump scheduling problem, sometimes referred to as the Optimal Water Flow (OWF) problem.
The OWF aims at determining a time-based policy for controllable components in a WDN, including pumps and valves, that minimizes the total cost of electrical power consumption.
From a mathematical programming perspective, the OWF is highly challenging for several reasons.
First, the OWF is constrained by energy loss and energy gain equations for pipes and pumps, respectively, which act as nonconvex constraints.
Second, the modeling of pumps and valves requires binary variables for indicating the statuses of these components, which increases the problem's combinatorial complexity.
Third, the modeling of time-evolving storage elements (e.g., elevated tanks) establishes a temporal dimension to the problem and increases its size considerably.

To alleviate these difficulties, this paper develops polyhedral relaxations of the OWF and a number of novel cutting plane and preprocessing methods that improve the primal-dual bound convergence of solution techniques.
These enhancements yield high-quality solutions that remain feasible to the original MINLP.
The contributions of this paper are as follows:
\begin{itemize}
    \item Tight polyhedral relaxation-based formulations of the canonical OWF;
    \item New valid inequalities inspired by duality-based constraint reformulations;
    \item Efficient relaxation-based optimization-based bound tightening (OBBT) algorithms;
    \item A novel optimization-based cut generation (OBCG) preprocessing algorithm;
    \item Empirical evaluations of formulation enhancements over a variety of OWF instances.
\end{itemize}
These contributions address two deficiencies in the current state of the art: (i) the slow improvement of dual bounds and (ii) the difficulty of generating feasible primal solutions.

The remainder of this paper proceeds as follows:
Section \ref{section:literature_review} reviews previous techniques from the WDN optimization literature;
Section \ref{section:problem_formulation} formulates the OWF as a MINLP;
Section \ref{section:convex_relaxations} develops polyhedral relaxations of the OWF;
Section \ref{section:valid-inequalities} states valid inequalities for the OWF, one of which is theoretically novel;
Section \ref{section:strengthening} develops a number of preprocessing algorithms that strengthen polyhedral relaxations;
Section \ref{section:computational_experiments} empirically evaluates our modeling, formulation, and algorithmic improvements;
and Section \ref{section:conclusion} concludes the paper.

\section{Literature Review}
\label{section:literature_review}
The study of WDN optimization resides within the broader category of optimization of \emph{nonlinear networks}, in which (i) flow is driven by potentials and (ii) potential loss along an edge is a nonlinear function of flow.
Critical infrastructure networks that satisfy these properties include potable water \citep{mala2017lost,mala2018lost}, natural gas \citep{rios2015optimization}, and crude oil \citep{sahebi2014strategic}, each of which has evolved as an important subfield within applied optimization.
\cite{raghunathan2013global} provides an overview of methods that have emerged to solve similar problems over the past five decades.
For the remainder of this paper, we reserve our discussion to WDNs, although we remark that the methods we describe could easily be extended to other nonlinear networks.

A thorough review of WDN operational optimization is provided by \cite{mala2017lost}, in which they note that the majority of studies consider the control of pumps and valves.
They also observe that most studies apply ``stochastic'' (e.g., evolutionary) heuristic techniques, while fewer use ``deterministic'' mathematical programming algorithms.
We remark that stochastic methods provide no guarantees of optimality nor, often times, feasibility.
Solutions obtained via such techniques could thus have deleterious effects if implemented in practice (e.g., permissiveness of pressure violations).
Finally, \cite{mala2017lost} find that studies typically consider WDNs with small numbers of nodes (or junctions), with $80\%$ concerning WDNs of one hundred nodes or fewer.
In line with the state of practice, we develop deterministic methods for similarly-sized WDNs.

More recently, mathematical programming techniques for optimal WDN operation have seen a resurgence in the literature.
For example, \cite{fooladivanda2017energy} consider a variant of the OWF similar to the one considered in this paper.
To solve the problem, they develop a mixed-integer second-order cone (SOC) relaxation and propose a specialized algorithm for a problem subclass.
In \cite{singh2019optimal}, a penalty-based SOC relaxation is developed.
In \cite{bonvin2017convex}, a heuristic using a similar relaxation strategy is developed.
Although these and other studies often assume WDNs with particular topologies (e.g., WDNs without loops) to provide feasibility guarantees, they have established the importance of \emph{convex relaxation} of nonlinearities that complicate the OWF problem.

A number of recent studies have established more general methods to solve the OWF and provide solutions closer to feasibility, independent of topology.
One example is given by \cite{vieira2020optimizing}, who develop a piecewise-linear relaxation of complicating nonlinearities, which is similarly done in this paper.
They also introduce network simplifications and valid inequalities, as well as an algorithm that attempts to recover feasible solutions from relaxation solutions using the de facto standard for WDN analysis, \textsc{EPANET}.
Another example is given by \cite{liu2020optimization}, who develop a mixed-integer linear programming (MILP) formulation for a WDN demand response application and consider a heuristic bound-tightening procedure using randomized network analysis results as bound proxies.

\cite{bonvin2021pump}, using an approach similar to ours, develop a MILP outer approximation of the OWF, as well as a linear programming/nonlinear programming branch and bound (LP/NLP-BB) algorithm that checks MINLP feasibility at integer-feasible nodes, similar to the algorithms used for WDN design by \cite{raghunathan2013global} and \cite{tasseff2020exact}.
They also employ a bound-tightening procedure to improve the strength of their MILP formulation.
Finally, they develop a heuristic that repairs physically infeasible solutions by permitting \emph{continuous}-duration pump activations.
However, because these heuristic solutions violate integrality assumptions, they are not feasible to the original OWF.

We consider \cite{bonvin2021pump} to be the current state of the art for OWF optimization.
However, differing from \cite{bonvin2021pump}, we develop outer \emph{and piecewise} MILP relaxations, which more accurately approximate nonlinear equations that model physical constraints.
Additionally, we evaluate several polyhedral relaxation-based OBBT techniques, whereas \cite{bonvin2021pump} employ a more computationally challenging MINLP-based OBBT.
Furthermore, this paper derives, computes, and applies several novel valid inequalities for the OWF.
Finally, this paper enforces integrality on pump and valve activations, which is a common feature of the canonical OWF described throughout the literature.

\section{Problem Formulation}
\label{section:problem_formulation}
Operation of a WDN involves the use of pumps to increase pressure and deliver water according to demands over a short-term planning horizon (e.g., twenty-four hours).
Together, pumps and tanks can be temporally coordinated to meet variable water demands while strategically achieving a desired objective, such as minimizing the total cost of energy consumed by the WDN.
An optimal ``pump schedule'' will often leverage the temporal variation in energy price over the planning horizon.
In this section, we formulate the OWF exactly as a MINLP.
Later, in Section \ref{section:convex_relaxations}, we present the proposed polyhedral relaxations.

\subsection{Water Network Modeling}
\label{section:network}
\subsubsection{Notation for Sets}
A WDN is represented by a directed graph $\mathcal{G} := (\mathcal{N}, \mathcal{L})$, where $\mathcal{N}$ is the set of nodes and $\mathcal{L}$ is the set of node-connecting components (links or arcs).
Since the OWF considers operational decisions across a fixed time horizon, $\mathcal{K} = \{1, 2, \dots, K\}$ is used to denote the set of time indices that describe the network's evolution.
The set $\tilde{\mathcal{K}} = \{1, 2, \dots, K + 1\}$ includes an additional time index, $K + 1$, to model final storage quantities.
We further define the sets of demands $\mathcal{D}$, reservoirs $\mathcal{R}$, and tanks $\mathcal{T}$ as disjoint subsets of the nodes, $\mathcal{N}$, in the network, where $\mathcal{D} \cup \mathcal{R} \cup \mathcal{T} = \mathcal{N}$, and where pass-through junctions with zero demand are also considered to reside in $\mathcal{D}$.
The set of node-connecting components $\mathcal{L}$ in the network comprises the disjoint subsets of pipes $\mathcal{A} \subset \mathcal{L}$, valves $\mathcal{V} \subset \mathcal{L}$, and fixed-speed pumps $\mathcal{P} \subset \mathcal{L}$, where $\mathcal{A} \cup \mathcal{V} \cup \mathcal{P} = \mathcal{L}$.
For convenience, we define the set of node-connecting components that are incident to node $i \in \mathcal{N}$, where $i$ is the tail (respectively, head) of the component, as being denoted by $\delta^{+}_{i} := \{(i, j) \in \mathcal{L}\}$ (respectively, $\delta^{-}_{i} := \{(j, i) \in \mathcal{L}\}$).

Next, we examine each of the component types independently, define their corresponding decision variables, and present constraints that each component enforces on the WDN's operations.
Specifically, for each component, we typically present two types of constraints: (i) operational limits and (ii) physical constraints.
We begin by examining the nodes, $\mathcal{N}$.

\subsubsection{Nodes}
Nodal potentials are denoted by the variables $h_{i}^{k}$, $i \in \mathcal{N}$, $k \in \tilde{\mathcal{K}}$, where each $h_{i}^{k}$ represents the total hydraulic head in units of length.
This quantity (hereafter referred to as ``head'') assimilates elevation and pressure heads.
Each head is constrained between lower and upper bounds, $\underline{h}_{i}^{k}$ and $\overline{h}_{i}^{k}$, respectively.
This implies the head constraints
\begin{equation}
    \underline{h}_{i}^{k} \leq h_{i}^{k} \leq \overline{h}_{i}^{k}, \, \forall i \in \mathcal{N}, \, \forall k \in \tilde{\mathcal{K}} \label{equation:node-head-bounds}.
\end{equation}

\paragraph{\textbf{Demands}}
Demands are nodes where water is supplied to end consumers.
Each demand is associated with a constant, $\overline{q}^{k}_{i}$, that denotes the demanded volumetric flow rate, where a negative value indicates consumption.
Without loss of generality, we also allow positive flows to represent injections (e.g., from a well).
The variables $q_{i}^{k} \in \mathbb{R}$, $i \in \mathcal{D}$, $k \in \mathcal{K}$, denote the flow consumed or supplied by each demand node.
When a demand is constant, $q_{i}^{k} = \overline{q}^{k}_{i}$.

\paragraph{\textbf{Reservoirs}}
Reservoirs are nodes where water is supplied to the WDN.
Each reservoir is modeled as an infinite source of flow with zero pressure head and constant elevation over a time step (i.e., $\underline{h}_{i}^{k} = \overline{h}_{i}^{k}$ at every reservoir, $i \in \mathcal{R}$).
Furthermore, the variables $q_{i}^{k} \geq 0$, $i \in \mathcal{R}$, $k \in \mathcal{K}$, are used to denote the \emph{outflow} of water from a reservoir $i$ at time index $k \in \mathcal{K}$.

\paragraph{\textbf{Tanks}}
Tanks store and discharge water over time.
Here, all tanks are assumed to be cylindrical with a fixed diameter $D_{i}$, $i \in \mathcal{T}$.
We assume there is no pressure head in the tanks, i.e., they are vented to the atmosphere.
The bottom of each tank, $B_{i}$, is located at or below the minimum water elevation, i.e., $B_{i} \leq \underline{h}_{i}^{k}$, $i \in \mathcal{T}$, and the maximum elevation of water is assumed to be $\overline{h}_{i}^{k}$.
The bounded variables $q_{i}^{k}$, $i \in \mathcal{T}$, $k \in \mathcal{K}$, denote the outflow (positive) or inflow (negative) through each tank.
The water volumes within the tanks are
\begin{equation}
    v_{i}^{k} := \frac{\pi}{4} D_{i}^{2} (h_{i}^{k} - B_{i}^{k}), \, \forall i \in \mathcal{T}, \, \forall k \in \tilde{\mathcal{K}} \label{equation:tank-volume-expression}.
\end{equation}
The Euler steps for integrating all tank volumes across time indices are then imposed with
\begin{equation}
    v_{i}^{k+1} = v_{i}^{k} - \Delta t^{k} q_{i}^{k}, \, \forall i \in \mathcal{T}, \, \forall k \in \mathcal{K} \label{equation:tank-volume-integration},
\end{equation}
where $\Delta t^{k}$ is the length of the time interval that connects times $k \in \mathcal{K}$ and $k + 1 \in \tilde{\mathcal{K}}$.

\subsubsection{Node-connecting Components}
\label{subsubsection:node-connecting-components}
Every link component $(i, j) \in \mathcal{L}$ is associated with a variable, $q_{ij}^{k}$, which denotes the volumetric flow rate across that component.
Assuming lower and upper bounds of $\underline{q}_{ij}^{k}$ and $\overline{q}_{ij}^{k}$, respectively, these variables are bounded via
\begin{equation}
    \underline{q}_{ij}^{k} \leq q_{ij}^{k} \leq \overline{q}_{ij}^{k}, \, \forall (i, j) \in \mathcal{L}, \, \forall k \in \mathcal{K} \label{equation:mincp-flow-bounds}.
\end{equation}
When $q_{ij}^{k}$ is positive, flow on $(i, j)$ is transported from node $i$ to $j$.
When $q_{ij}^{k}$ is negative, flow is transported from node $j$ to $i$.
At zero flow, the flow direction is considered ambiguous.

\paragraph{\textbf{Pipes}}
Pipes are the primary means for transporting water in a WDN.
Water flowing through a pipe will exhibit frictional loss due to contact with the pipe wall.
In this paper, energy loss relationships that link pipe flow and head (i.e., ``head loss'' equations) are modeled by the Hazen-Williams equation \citep{ormsbee2016darcy}, requiring the constraints
\begin{equation}
    h_{i}^{k} - h_{j}^{k} = L_{ij} r_{ij} q_{ij}^{k} \left\lvert q_{ij}^{k} \right\rvert^{\alpha - 1}, \, \forall (i, j) \in \mathcal{A}, \, \forall k \in \mathcal{K} \label{equation:mincp-pipe-head-loss}.
\end{equation}
Here, $\alpha = 1.852$ for each relationship and $r_{ij}$ denotes the resistance per unit length, which comprises all length-independent constants that appear in the Hazen-Williams equation.

\paragraph{\textbf{Valves}}
Valves control the flow of water to specific portions of the WDN.
Here, valves are elements that are either open or closed.
The operating status of each valve $(i, j) \in \mathcal{V}$ is indicated using a binary variable, $z_{ij}^{k} \in \{0, 1\}$, where $z_{ij}^{k} = 1$ corresponds to an open valve and $z_{ij}^{k} = 0$ to a closed valve.
These binary variables restrict the flow across each valve as
\begin{equation}
    \underline{q}_{ij}^{k} z_{ij}^{k} \leq q_{ij}^{k} \leq \overline{q}_{ij}^{k} z_{ij}^{k}, \, z_{ij}^{k} \in \{0, 1\}, \, \forall (i, j) \in \mathcal{V}, \, \forall k \in \mathcal{K} \label{equation:mincp-valve-flow-bounds}.
\end{equation}
Furthermore, when a valve is open, the heads at the nodes connected by that valve are equal.
When the valve is closed, these heads are decoupled.
This disjunctive phenomenon is modeled via the following set of constraints involving the binary indicator variables $z_{ij}^{k}$:
\begin{equation}
    (1 - z_{ij}^{k}) (\underline{h}_{i}^{k} - \overline{h}_{j}^{k}) \leq h_{i}^{k} - h_{j}^{k} \leq (1 - z_{ij}^{k}) (\overline{h}_{i}^{k} - \underline{h}_{j}^{k}), \, \forall (i, j) \in \mathcal{V}, \, \forall k \in\mathcal{K} \label{equation:mincp-valve-head}.
\end{equation}
That is, if $z_{ij}^{k} = 1$, then $h_{i}^{k} = h_{j}^{k}$.
Otherwise, if $z_{ij}^{k} = 0$, then $h_{i}^{k}$ and $h_{j}^{k}$ are decoupled.

\paragraph{\textbf{Pumps}}
Each pump $(i, j) \in \mathcal{P}$ increases the head from node $i$ to $j$ when active and permits only unidirectional flow.
Here, we consider fixed-speed pumps.
When the pump is off, there is zero flow, and heads at adjacent nodes are decoupled.
When the pump is on, there is positive flow (greater than or equal to some fixed $\underline{q}_{ij}^{k +}$),
and the head increase from $i$ to $j$ is modeled by a nonlinear function.
The variable $z_{ij}^{k} \in \{0, 1\}$ indicates the status of each pump, where $z_{ij}^{k} = 1$ if $q_{ij}^{k} \geq \underline{q}_{ij}^{k +}$ and $z_{ij}^{k} = 0$ if $q_{ij}^{k} < \underline{q}_{ij}^{k +}$.
This implies the disjunctive bounds
\begin{equation}
    \underline{q}_{ij}^{k} = 0 \leq \underline{q}_{ij}^{k +} z_{ij}^{k} \leq q_{ij}^{k} \leq \overline{q}_{ij}^{k} z_{ij}^{k}, \, z_{ij}^{k} \in \{0, 1\}, \, \forall (i, j) \in \mathcal{P}, \, \forall k \in \mathcal{K} \label{equation:mincp-pump-flow-bounds}.
\end{equation}
The variable $g_{ij}^{k} \geq 0$ is introduced for each pump to denote the head increase (or gain) that results from that pump. 
Similar to the standard for water network analysis, \textsc{EPANET} \citep{Ros2000}, we model pump head gains via strictly concave equations of the form
\begin{equation}
    a_{ij} z_{ij}^{k} + b_{ij} (q_{ij}^{k})^{c_{ij}} = g_{ij}^{k}, \, \forall (i, j) \in \mathcal{P}, \, \forall k \in \mathcal{K} \label{equation:mincp-pump-head-gain}.
\end{equation}
In these constraints, $a_{ij}$ and $c_{ij}$ are positive constants, and $b_{ij}$ is a negative constant.
Further, note that $q_{ij}^{k}$ and $a_{ij} z_{ij}^{k}$ are restricted to zero when $z_{ij}^{k}$ is zero.
This ensures that when a pump is off, the corresponding head gain $g_{ij}^{k}$ is also zero.
To ensure the decoupling of hydraulic heads when a pump is off, the following disjunctive constraints are employed:
\begin{subequations}
\begin{align}
    h_{i}^{k} - h_{j}^{k} + g_{ij}^{k} &\leq (1 - z_{ij}^{k}) (\overline{h}_{i}^{k} - \underline{h}_{j}^{k}), \, \forall (i, j) \in \mathcal{P}, \, \forall k \in \mathcal{K} \\
    h_{i}^{k} - h_{j}^{k} + g_{ij}^{k} &\geq (1 - z_{ij}^{k}) (\underline{h}_{i}^{k} - \overline{h}_{j}^{k}), \, \forall (i, j) \in \mathcal{P}, \, \forall k \in \mathcal{K}.
\end{align}
\label{equation:mincp-pump-head}%
\end{subequations}
Note that when $z_{ij}^{k} = 1$, the pump is on, and the head \emph{gain} between the two nodes is $g_{ij}^{k}$.

\subsubsection{Flow Conservation}
Satisfaction of flow (or mass) conservation requires flow balance constraints to be enforced at every node $i \in \mathcal{N}$ across the time indices $\mathcal{K}$.
That is,
\begin{equation}
    \sum_{\mathclap{(i, j) \in \delta^{+}_{i}}} q_{ij}^{k} - \sum_{\mathclap{(j, i) \in \delta^{-}_{i}}} q_{ji}^{k} = q_{i}^{k}, \, \forall i \in \mathcal{N}, \, \forall k \in \mathcal{K}.
\label{equation:flow-conservation}%
\end{equation}

\subsection{Optimal Water Flow}
\label{section:optimal-water-flow}
Section \ref{section:network} described the variables and constraints required to model the physics and feasible operation of a WDN.
However, a feasible operational schedule is not necessarily optimal, and there is a combinatorial number of possibilities for pump and valve settings, $z_{ij}^{k}$, at each time index.
These are the primary OWF decision variables of interest.
As described in Section \ref{section:introduction}, an optimal operational schedule will minimize the cost of energy consumption of pumps in the network.
A solution to the OWF should also satisfy a number of other common, practical criteria.
In the following subsections, we enumerate these criteria.

\subsubsection{Tank Volume Recovery}
At the end of the fixed operational schedule, the volume of water within each tank should be at least as large as its initial volume.
This ensures that (i) the operational schedule can be repeatedly applied to subsequent planning periods with similar demand profiles and (ii) tank volumes at the very end of the planning period are not so low that operation in subsequent periods cannot be as efficient.
This implies
\begin{equation}
    v_{i}^{1} \leq v_{i}^{K + 1}, \, \forall i \in \mathcal{T} \label{equation:tank-volume-recovery},
\end{equation}
where, in most of the literature, as well as in this paper, $v_{i}^{1}$ is assumed to be a known, fixed value.
Denoting each fixed value by $\tilde{v}_{i}$, this would further imply that $v_{i}^{1} = \tilde{v}_{i}$ for $i \in \mathcal{T}$.

\subsubsection{Pump Switching Limits}
Maintenance accounts for roughly $10\%$ of the net present value lifecycle cost of a pump \citep{nault2015lifecycle}.
Frequent activation and deactivation (i.e., switching) generally reduces the pump's lifetime and increases its overall maintenance cost.
Thus, WDN operators can require operational constraints to limit the number of pump activations and deactivations.
Assuming that each pump in the network must \emph{remain on} for a minimum time of $\tau^{\textrm{on}}$, \emph{off} for a minimum time of $\tau^{\textrm{off}}$, and \emph{switched on} no more than $N_{ij}$ times, we adopt the following constraints of \citet{ghaddar2015lagrangian}:
\begin{subequations}
\begin{align}
    z_{ij}^{\textrm{on}, k}, z_{ij}^{\textrm{off}, k} \in \{0, 1\}, \, &\forall (i, j) \in \mathcal{P}, \, \forall k \in \mathcal{K} \\
    z_{ij}^{k} - z_{ij}^{k - 1} \leq z_{ij}^{\textrm{on}, k}, \, &\forall (i, j) \in \mathcal{P}, \, \forall k \in \mathcal{K} \setminus \{1\} \label{equation:pump-switch-1} \\
    z_{ij}^{\textrm{on}, k} \leq z_{ij}^{k^{\prime}}, \, &\forall (i, j) \in \mathcal{P}, \, \forall k \in \mathcal{K}, \, \forall k^{\prime} \in t(k) \leq t(k^{\prime}) \leq t(k) + \tau^{\textrm{on}} \label{equation:pump-switch-2} \\
    z_{ij}^{k - 1} - z_{ij}^{k} \leq z_{ij}^{\textrm{off}, k}, \, &\forall (i, j) \in \mathcal{P}, \, \forall k \in \mathcal{K} \setminus \{1\} \label{equation:pump-switch-3} \\
    z_{ij}^{k^{\prime}} \leq 1 - z_{ij}^{\textrm{off}, k}, \, &\forall (i, j) \in \mathcal{P}, \, \forall k \in \mathcal{K}, \, \forall k^{\prime} \in t(k) \leq t(k^{\prime}) \leq t(k) + \tau^{\textrm{off}} \label{equation:pump-switch-4} \\
    \sum_{k \in \mathcal{K}} z_{ij}^{\textrm{on}, k} \leq N_{ij}, \, &\forall (i, j) \in \mathcal{P} \label{equation:pump-switch-5}.
\end{align}
\label{equation:pump-switch}%
\end{subequations}
Here, $t(k)$ returns the time at index $k \in \mathcal{K}$, while $z_{ij}^{\textrm{on}, k} \in \{0, 1\}$ and $z_{ij}^{\textrm{off}, k} \in \{0, 1\}$ indicate whether pump $(i, j) \in \mathcal{P}$ has been switched \emph{on} or \emph{off} at $k$, respectively.
Constraints \eqref{equation:pump-switch-1} ensure that $z_{ij}^{\textrm{on}, k}$ is equal to one when a pump \emph{becomes} active at $k$.
Constraints \eqref{equation:pump-switch-2} ensure that, if a pump has been switched on, it remains on for at least a duration of $\tau^{\textrm{on}}$.
Constraints \eqref{equation:pump-switch-3} ensure that $z_{ij}^{\textrm{off}, k}$ is equal to one when a pump \emph{becomes inactive} at $k$.
Constraints \eqref{equation:pump-switch-4} ensure that, if a pump has been switched off, it remains off for at least a duration of $\tau^{\textrm{off}}$.
Finally, Constraints \eqref{equation:pump-switch-5} limit the total number of \emph{on} switches to $N_{ij}$.

\subsubsection{Cost Minimization}
The final criterion for an OWF solution is that the operational decisions minimize the total cost of power consumption over the planning horizon.
This cost arises from the energy requirements of pumps.
For the considered fixed-speed pumps, we utilize a power curve that is a linear function of flow.
This approximation has been used in the literature (e.g., \citealp{bonvin2021pump}) and further implies the cost function
\begin{equation}
    f(q, z) = \sum_{k \in \mathcal{K}} \sum_{(i, j) \in \mathcal{P}} \lambda_{ij}^{k} q_{ij}^{k} + \mu_{ij}^{k} z_{ij}^{k} \label{equation:owf-objective-function},
\end{equation}
where $\lambda_{ij}^{k}$ and $\mu_{ij}^{k}$ are constant parameters that incorporate both the modeling of power consumption at pump $(i, j) \in \mathcal{P}$ as well as the price of electricity at time index $k \in \mathcal{K}$.

\subsubsection{The Optimal Water Flow Problem}
Combining feasibility, tank recovery, pump switching limits, and minimization of cost gives rise to the canonical OWF problem,
\begin{equation}
\tag{MINLP}
\begin{aligned}
    & \text{minimize}
    & & \textnormal{Cost of electricity} && f(q, z) ~ \textnormal{of Equation} ~ \eqref{equation:owf-objective-function} \\
    & \text{subject to}
    & & \textnormal{Operational feasibility} && \textnormal{Constraints} ~ \eqref{equation:node-head-bounds}, \eqref{equation:tank-volume-integration}{\textnormal{--}}\eqref{equation:flow-conservation} \\
    & & & \textnormal{Tank volume recovery} && \textnormal{Constraints} ~ \eqref{equation:tank-volume-recovery} \\
    & & & \textnormal{Pump switching limits} && \textnormal{Constraints} ~ \eqref{equation:pump-switch}.
\end{aligned}
\label{equation:mincp}%
\end{equation}

\section{Polyhedral Relaxations}
\label{section:convex_relaxations}
There are two sources of nonlinear nonconvexity that render \eqref{equation:mincp} difficult to solve directly.
These are Constraints \eqref{equation:mincp-pipe-head-loss} and \eqref{equation:mincp-pump-head-gain}, which model each pipe's head loss and each pump's head gain, respectively.
Numerical methods capable of solving even \emph{convex} MINLPs to global optimality, e.g., \textsc{Juniper} \citep{juniper} and \textsc{BONMIN} \citep{bonami2007bonmin}, are not guaranteed to do so efficiently.
Fractional exponents in head loss and gain relationships exacerbate this issue.
A number of recent studies have addressed these challenges via polyhedral (MILP) relaxations of nonlinear constraints \citep{vieira2020optimizing,bonvin2021pump,liu2020optimization}.
Unlike MINLP-based formulations, these relaxations can more easily exploit the efficiency of commercial MILP solvers.
In the following subsections, we extend prior methods to develop two new flow direction-based MILP relaxations that employ novel flow bounds based on the activation and flow direction of node-connecting components.
These modifications enable the further strengthening of relaxed formulations.

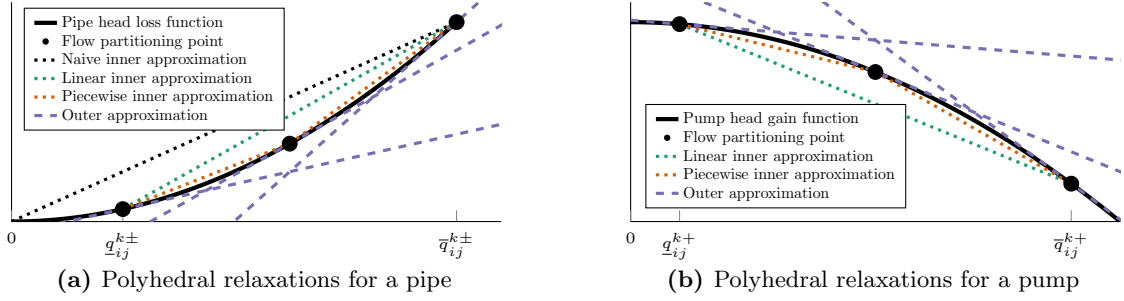
\begin{figure}[t]
    \centering
    \begin{subfigure}[t]{0.49\linewidth}
        \centering
        \begin{tikzpicture}[declare function={f_1(\x)=\x^(2);}]
	\begin{axis}[axis lines*=middle,xtick={0.25,1.0},height=4.5cm,width=\linewidth,ticklabel style = {font=\tiny},
	             legend cell align={left},legend style={nodes={scale=0.5,transform shape}},ymin=0.0,
	             xmin=0.0,xmax=1.1,ymax=1.1,extra x ticks={0},
	             legend style={at={(0.025, 0.975)},anchor=north west},
	             legend image post style={scale=0.5},ymajorticks=false,
	             xticklabels={$\underline{q}_{ij}^{k \pm}$,$\overline{q}_{ij}^{k \pm}$}]
    \addplot[ultra thick,black,samples=50,domain=0.0:1.0,name path=f] {f_1(x)};
    \addlegendentry{Pipe head loss function};
    
    \addlegendimage{only marks,mark=*,black};
    \addlegendentry{Flow partitioning point};
    
    \draw[name path=inner_naive,very thick,dotted,color=black] (axis cs:0.0,0.0) -- (axis cs:1.0,1.0);
    \addlegendimage{line width=0.3mm,very thick,dotted,color=black};
    \addlegendentry{Naive inner approximation};
    
    \draw[name path=inner_standard,very thick,dotted,index of colormap=0 of Dark2-8] (0.25,0.0625) -- (1.0,1.0);
    \addlegendimage{line width=0.3mm,very thick,dotted,index of colormap=0 of Dark2-8};
    \addlegendentry{Linear inner approximation};
    
    \draw[name path=inner_pw_1,very thick,dotted,index of colormap=1 of Dark2-8] (0.25,0.0625) -- (0.625,0.390625);
    \draw[name path=inner_pw_2,very thick,dotted,index of colormap=1 of Dark2-8] (0.625,0.390625) -- (1.0,1.0);
    \addlegendimage{line width=0.3mm,very thick,dotted,index of colormap=1 of Dark2-8};
    \addlegendentry{Piecewise inner approximation};
    
    \draw[name path=outer_1,very thick,dashed,index of colormap=2 of Dark2-8] (axis cs:0.0,-0.0625) -- (axis cs: 1.25,0.5625);
    \draw[name path=outer_1,very thick,dashed,index of colormap=2 of Dark2-8] (axis cs:0.0,-0.390625) -- (axis cs: 1.25,1.171875);
    \draw[name path=outer_1,very thick,dashed,index of colormap=2 of Dark2-8] (axis cs:0.0,-1) -- (axis cs: 1.25,1.5);
    \addlegendimage{line width=0.3mm,very thick,dashed,index of colormap=2 of Dark2-8};
    \addlegendentry{Outer approximation};

    \fill (0.25,0.0625)    circle [radius=3.0pt];
    \fill (0.625,0.390625) circle [radius=3.0pt];
    \fill (1.0,1.0)        circle [radius=3.0pt];
	\end{axis}
\end{tikzpicture}
        \caption{Polyhedral relaxations for a pipe}
        \label{subfigure:polyhedral-pipe}
    \end{subfigure}
    \begin{subfigure}[t]{0.49\linewidth}
        \centering
        \begin{tikzpicture}[declare function={f_1(\x)=-\x^(2) + 1;}]
	\begin{axis}[axis lines*=middle,xtick={0.10,0.90},height=4.5cm,width=\linewidth,
	             legend cell align={left},legend style={nodes={scale=0.5,transform shape}},ticklabel style = {font=\tiny},
	             xmin=0.0,xmax=1.0,ymin=0.0,ymax=1.1,extra x ticks={0},
	             legend style={at={(0.03, 0.06)},anchor=south west},
	             legend image post style={scale=0.5},ymajorticks=false,
	             xticklabels={$\underline{q}_{ij}^{k +}$,$\overline{q}_{ij}^{k +}$}]
    \addplot[ultra thick,black,samples=50,domain=0.0:1.0,name path=f] {f_1(x)};
    \addlegendentry{Pump head gain function};
    
    \addlegendimage{only marks,mark=*,black};
    \addlegendentry{Flow partitioning point};
    
    
    \draw[name path=inner_standard,very thick,dotted,index of colormap=0 of Dark2-8] (0.1,0.99) -- (0.9,0.19);
    \addlegendimage{line width=0.3mm,very thick,dotted,index of colormap=0 of Dark2-8};
    \addlegendentry{Linear inner approximation};
    
    \draw[name path=inner_pw_1,very thick,dotted,index of colormap=1 of Dark2-8] (0.1,0.99) -- (0.5,0.75);
    \draw[name path=inner_pw_2,very thick,dotted,index of colormap=1 of Dark2-8] (0.5,0.75) -- (0.9,0.19);
    \addlegendimage{line width=0.3mm,very thick,dotted,index of colormap=1 of Dark2-8};
    \addlegendentry{Piecewise inner approximation};
    
    \draw[name path=outer_1,very thick,dashed,index of colormap=2 of Dark2-8] (axis cs:0.0,1.01) -- (axis cs: 1.25,0.76);
    \draw[name path=outer_1,very thick,dashed,index of colormap=2 of Dark2-8] (axis cs:0.0,1.25) -- (axis cs: 1.25,0.0);
    \draw[name path=outer_1,very thick,dashed,index of colormap=2 of Dark2-8] (axis cs:0.0,1.81) -- (axis cs: 1.25,-0.44);
    \addlegendimage{line width=0.3mm,very thick,dashed,index of colormap=2 of Dark2-8};
    \addlegendentry{Outer approximation};

    \fill (0.1,0.99) circle [radius=3.0pt];
    \fill (0.5,0.75) circle [radius=3.0pt];
    \fill (0.9,0.19) circle [radius=3.0pt];
	\end{axis}
\end{tikzpicture}
        \caption{Polyhedral relaxations for a pump}
        \label{subfigure:polyhedral-pump}
    \end{subfigure}
    \vspace{0.5em}
    \caption{Illustrations of polyhedral relaxations for the pipe/pump head loss/gain constraints considered in this paper. Dotted lines indicate separate possible \emph{inner approximations} of the nonlinear functions, and dashed lines indicate \emph{outer approximations}. Points indicate function partitioning points selected a priori.}
    \label{figure:polyhedral_relaxations}
\end{figure}

Figure \ref{figure:polyhedral_relaxations} illustrates the direction-based relaxations formulated in the following subsections.
Here, solid lines correspond to the original pipe head loss and pump head gain constraints, decomposed by flow direction (i.e., $+$ or $-$).
Both formulations are based on linear \emph{outer approximations} of the complicating nonlinearities, indicated by dashed purple lines.
A number of possible \emph{inner approximations}, denoted by dotted black, green, and orange lines, can bound these outer approximations.
Note that the novelty of directed lower bounds, i.e., $\underline{q}_{ij}^{\pm}$, helps reduce the size of the feasible regions of these relaxations.
Second, note that piecewise inner approximations, illustrated with dotted orange lines, can substantially strengthen relaxations, albeit at a greater computational cost.
In the following subsections, we algebraically formulate the polyhedral relaxations illustrated in Figure \ref{figure:polyhedral_relaxations}.

\subsection{Flow Direction-based Decomposition}
\label{section:micp-relaxation}
We begin by introducing the variables $q_{ij}^{k \pm}$ to denote nonnegative flows in the two possible flow directions along each link, $(i, j) \in \mathcal{L}$.
This decomposition first implies the definition
\begin{equation}
    q_{ij}^{k} := q^{k+}_{ij} - q^{k-}_{ij}, \, \forall (i, j) \in \mathcal{L}, \, \forall k \in \mathcal{K} \label{equation:flow-equations}.
\end{equation}
Next, discrete variables $y_{ij}^{k} \in \{0, 1\}$, $(i, j) \in \mathcal{L}$, are introduced to model the flow direction, where $y_{ij}^{k} = 1$ implies flow from $i$ to $j$ and $y_{ij}^{k} = 0$ implies flow from $j$ to $i$.
It follows that
\begin{subequations}
\begin{align}
    y_{ij}^{k} \in \{0, 1\}, \; &\forall (i, j) \in \mathcal{L}, \, \forall k \in \mathcal{K} \\
    \underline{q}_{ij}^{k+} y_{ij}^{k} \leq q_{ij}^{k+} \leq \overline{q}_{ij}^{k+} y_{ij}^{k}, \; &\forall (i, j) \in \mathcal{L}, \, \forall k \in \mathcal{K} \\
    \underline{q}_{ij}^{k-} (1 - y_{ij}^{k}) \leq q_{ij}^{k-} \leq \overline{q}_{ij}^{k-} (1 - y_{ij}^{k}), \; &\forall (i, j) \in \mathcal{L}, \, \forall k \in \mathcal{K},
\end{align}
\label{equation:micpr-direction-bounds}%
\end{subequations}
where $\underline{q}_{ij}^{k \pm} \geq 0$, $\overline{q}_{ij}^{k +} = \max\{0, \overline{q}_{ij}^{k}\}$, and $\overline{q}_{ij}^{k -} = \max\{0, -\underline{q}_{ij}^{k}\}$.
Compared to recent relaxations, a subtle but important difference is the modeling of \emph{directed lower flow bounds} $\underline{q}_{ij}^{k \pm}$, which can be conservatively estimated as zero.
However, as will later be seen, methods can be used to increase these directed bounds, thus tightening direction-based relaxations.

For pipes, nonnegative head difference variables $\Delta h_{ij}^{k \pm}$ are used to model losses via
\begin{subequations}
\begin{align}
    \Delta h_{ij}^{k+} - \Delta h_{ij}^{k-} = h_{i}^{k} - h_{j}^{k}, \; &\forall (i, j) \in \mathcal{A}, \, \forall k \in \mathcal{K} \label{equation:micpr-head-equality} \\
    0 \leq \Delta h_{ij}^{k+} \leq y_{ij}^{k} (\max\{0, \overline{h}_{i}^{k} - \underline{h}_{j}^{k}\}), \; &\forall (i, j) \in \mathcal{A}, \, \forall k \in \mathcal{K} \\
    0 \leq \Delta h_{ij}^{k-} \leq (1 - y_{ij}^{k}) (\max\{0, \overline{h}_{j}^{k} - \underline{h}_{i}^{k}\}), \; &\forall (i, j) \in \mathcal{A}, \, \forall k \in \mathcal{K}.
\end{align}%
\label{equation:micpr-head-difference-bounds}%
\end{subequations}
We also employ valid inequalities to exploit knowledge concerning flow directions in the WDN.
We omit them here, although they are detailed in the supplemental material.

\subsection{Mixed-integer Linear Relaxation}
\label{section:mixed-integer-linear-relaxation}
In this subsection, we define the non-piecewise MILP relaxation of the OWF, where nonlinearities are bounded by the \emph{dashed purple} and \emph{dotted green} lines in Figure \ref{figure:polyhedral_relaxations}.
This relaxation \emph{outer-approximates} nonlinearities, a well-known technique dating back to \citet{10.2307/2099058}.

\subsubsection{Linearization of Pipe Head Loss Constraints}
To relax Constraints \eqref{equation:mincp-pipe-head-loss}, discrete sets of flow partitioning points, $\mathcal{Q}_{ij}^{k \pm}$, are first assumed to have been defined a priori for $(i, j) \in \mathcal{A}$, $k \in \mathcal{K}$.
Constraints \eqref{equation:mincp-pipe-head-loss} are then linearly outer-approximated via the constraints
\begin{equation}
    r_{ij} [(\hat{q}_{ij}^{k \pm})^{\alpha} y_{ij}^{k \pm} +
        \alpha (\hat{q}_{ij}^{k \pm})^{\alpha - 1}
        (q_{ij}^{k \pm} - \hat{q}_{ij}^{k \pm} y_{ij}^{k \pm})] \leq
        \frac{\Delta h_{ij}^{k \pm}}{L_{ij}}, \forall (i, j) \in \mathcal{A}, \,
        \forall k \in \mathcal{K}, \, \forall \hat{q}_{ij}^{k \pm}
        \in \mathcal{Q}_{ij}^{k \pm},
\label{equation:milp-head-loss}%
\end{equation}
where $y_{ij}^{k +} := y_{ij}^{k}$ and $y_{ij}^{k -} := 1 - y_{ij}^{k}$ are defined to formulate the constraints more concisely.
These constraints are illustrated by the dashed purple lines in Figure \ref{subfigure:polyhedral-pipe}.
Observing the lower boundedness of these head loss relaxations, we also introduce linear upper bounds on the constraints that employ limits of $q_{ij}^{k \pm}$.
These linear upper-bounding constraints are
\begin{equation}
     r_{ij} \left([\underline{q}_{ij}^{k \pm}]^{\alpha} y_{ij}^{k \pm} + \frac{[(\underline{q}_{ij}^{k \pm})^{\alpha} - (\overline{q}_{ij}^{k \pm})^{\alpha}] [q_{ij}^{k \pm} - \underline{q}_{ij}^{k \pm} y_{ij}^{k \pm}]}{\underline{q}_{ij}^{k \pm} - \overline{q}_{ij}^{k \pm}} \right) \geq \frac{\Delta h_{ij}^{k \pm}}{L_{ij}}, \, \forall (i, j) \in \mathcal{A}, \, \forall k \in \mathcal{K} \label{equation:micpr-head-loss-ub}.%
\end{equation}
These upper bounds are illustrated by the dotted green line in Figure \ref{subfigure:polyhedral-pipe}.
Note that Constraints \eqref{equation:milp-head-loss} and \eqref{equation:micpr-head-loss-ub} are \emph{activated} by the values of $y_{ij}^{k \pm}$.
Furthermore, the strength of each linear upper-bounding constraint is determined by the proximity of $\underline{q}_{ij}^{k \pm}$ and $\overline{q}_{ij}^{k \pm}$.

\subsubsection{Linearization of Pump Head Gain Constraints}
For pumps, discrete sets $\mathcal{Q}_{ij}^{k +}$ are similarly used to linearly outer-approximate the head gain Constraints \eqref{equation:mincp-pump-head-gain} via
\begin{equation}
    \begin{gathered}
        g_{ij}^{k} \leq [a_{ij} + b_{ij} (\hat{q}_{ij}^{k +})^{c_{ij}}] z_{ij}^{k} + b_{ij} c_{ij} (\hat{q}_{ij}^{k +})^{c_{ij} - 1} (q_{ij}^{k +} - \hat{q}_{ij}^{k +} z_{ij}^{k}), \\
        \forall (i, j) \in \mathcal{P}, ~ \forall k \in \mathcal{K}, ~ \forall \hat{q}_{ij}^{k +} \in \mathcal{Q}_{ij}^{k +} \label{equation:milp-head-gain}.
    \end{gathered}
\end{equation}
These constraints are illustrated by the dashed purple lines in Figure \ref{subfigure:polyhedral-pump}.
Recall that $\overline{q}_{ij}^{k-} = q_{ij}^{k-} = 0$ for all $(i, j) \in \mathcal{P}$, $k \in \mathcal{K}$, and that $g_{ij}^{k} = 0$ when $z_{ij}^{k}$ (and thus $q_{ij}^{k+}$) is equal to zero.
Similar to the \emph{linear upper-bounding} Constraints \eqref{equation:micpr-head-loss-ub} for head loss relaxations, the \emph{linear lower bounds} for the pump head gain relaxation Constraints \eqref{equation:milp-head-gain} are stated as
\begin{equation}
    g_{ij}^{k} \geq (a_{ij} + b_{ij} [\underline{q}_{ij}^{k +}]^{c_{ij}}) z_{ij}^{k} + \frac{b_{ij} ([\underline{q}_{ij}^{k +}]^{c_{ij}} - [\overline{q}_{ij}^{k +}]^{c_{ij}}) (q_{ij}^{k +} - \underline{q}_{ij}^{k +} z_{ij}^{k})}{\underline{q}_{ij}^{k +} - \overline{q}_{ij}^{k +}}, \, \forall (i, j) \in \mathcal{P}, \, \forall k \in \mathcal{K} \label{equation:micpr-pump-head-gain-lb}.
\end{equation}
These constraints are illustrated by the dotted green line in Figure \ref{subfigure:polyhedral-pump}.
Similar to Constraints \eqref{equation:micpr-head-loss-ub}, the strength of each Constraint \eqref{equation:micpr-pump-head-gain-lb} is dependent on pump flow bounds.

\subsubsection{Outer Approximation of the Optimal Water Flow Problem}
Given the previously-described constraints, the MILP outer approximation of the OWF is written as
\begin{equation}
\tag{MILP-OA}
\begin{aligned}
    & \text{minimize}
    & & \textnormal{Cost of electricity} && f(q, z) ~ \textnormal{of Equation} ~ \eqref{equation:owf-objective-function} \\
    & \text{subject to}
    & & \textnormal{\eqref{equation:mincp} subset} && \textnormal{Constraints} ~ \eqref{equation:node-head-bounds}, \eqref{equation:tank-volume-integration}, \eqref{equation:mincp-valve-flow-bounds}\textnormal{--}\eqref{equation:mincp-pump-flow-bounds}, \eqref{equation:mincp-pump-head}\textnormal{--}\eqref{equation:pump-switch} \\
    & & & \textnormal{Relaxed physics} && \textnormal{Constraints} ~ \eqref{equation:micpr-direction-bounds}\textnormal{--}\eqref{equation:micpr-pump-head-gain-lb}.
\end{aligned}
\label{equation:milp-oa}%
\end{equation}
where the definitions of $q_{ij}^{k}$, $(i, j) \in \mathcal{L}$, $k \in \mathcal{K}$, from Equations \eqref{equation:flow-equations} are used where necessary.

\subsection{Mixed-integer Piecewise Linear Relaxation}
\label{section:mixed-integer-piecewise-linear-relaxation}
To improve the relaxation strength of \eqref{equation:milp-oa}, this section develops a mixed-integer \emph{piecewise} linear relaxation of \eqref{equation:mincp} that leverages outer \emph{and} piecewise-inner approximations of nonlinear constraints.
Compared to \eqref{equation:milp-oa}, the approach is similar but forms tighter piecewise envelopes of nonlinear constraints rather than simpler linear bounds.
These differences are illustrated by comparing the piecewise inner-approximating \emph{orange lines} of Figure \ref{figure:polyhedral_relaxations} with the non-piecewise, linear inner-approximating green lines.

\subsubsection{Linearization of Pipe Head Loss Constraints}
\label{section:linearization-head-loss}%
To strengthen \eqref{equation:milp-oa}, we introduce piecewise \emph{inner} (or upper) approximations of head loss constraints, which are stronger than the upper-bounding Constraints \eqref{equation:micpr-head-loss-ub}.
Letting $0 \leq \lambda_{ijp}^{k \pm} \leq 1$, $p \in \{1, 2, \dots, \lvert \mathcal{Q}_{ij}^{k \pm} \rvert\}$, $(i, j) \in \mathcal{A}$, $k \in \mathcal{K}$, denote the continuous convex combination variables and $x_{ijp}^{k \pm} \in \{0, 1\}$, $p \in \{2, \dots, \lvert \mathcal{Q}_{ij}^{k \pm} \lvert\}$, denote the binary convex combination variables, these are
\begin{equation}
     \sum_{p = 1}^{\mathclap{\lvert \mathcal{Q}_{ij}^{k \pm} \rvert}} r_{ij} (\hat{q}_{ijp}^{k \pm})^{\alpha} \lambda_{ijp}^{k \pm} \geq \frac{\Delta h_{ij}^{k \pm}}{L_{ij}}, \, \forall (i, j) \in \mathcal{A}, \, \forall k \in \mathcal{K}.
\label{equation:milp-convex-combination-pipe-loss}%
\end{equation}
The directed flow variables are similarly constrained by the convex combination, i.e.,
\begin{equation}
    q_{ij}^{k \pm} = \sum_{p = 1}^{\mathclap{\lvert \mathcal{Q}_{ij}^{k \pm} \rvert}} \hat{q}_{ijp}^{k \pm} \lambda_{ijp}^{k \pm}, \, \forall (i, j) \in \mathcal{A}, \, \forall k \in \mathcal{K}.
\label{equation:milp-convex-combination-pipe-flow}%
\end{equation}
The activations of direction-dependent convex combination variables are then limited by
\begin{subequations}
\begin{align}
    \sum_{p = 1}^{\mathclap{\lvert \mathcal{Q}_{ij}^{k \pm} \rvert}} \lambda_{ijp}^{k \pm} = y_{ij}^{k \pm},
    ~ \sum_{p = 2}^{\mathclap{\lvert \mathcal{Q}_{ij}^{k \pm} \rvert}} x_{ijp}^{k \pm} = y_{ij}^{k \pm}, \; &\forall (i, j) \in \mathcal{A}, \, \forall k \in \mathcal{K}
\label{equation:milp-convex-combination-direction} \\
    \lambda_{ijp}^{k \pm} \leq x_{ijp}^{k \pm} + x_{ij, p + 1}^{k \pm}, \; &\forall (i, j) \in \mathcal{A}, \, \forall k \in \mathcal{K}, \, \forall p \in \{2, 3, \dots, \lvert \mathcal{Q}_{ij}^{k \pm} \rvert - 1\} \label{equation:milp-convex-combination-direction1} \\
    \lambda_{ij, 1}^{k \pm} \leq x_{ij, 2}^{k \pm}, \; \lambda_{ij, \lvert \mathcal{Q}_{ij}^{k \pm} \rvert}^{k \pm} \leq x_{ij, \lvert \mathcal{Q}_{ij}^{k \pm} \rvert}^{k \pm}, \; &\forall (i, j) \in \mathcal{A}, \, \forall k \in \mathcal{K}.\label{equation:milp-convex-combination-direction2}
\end{align}
\label{equation:milp-convex-combination-relations}%
\end{subequations}
where Constraints \eqref{equation:milp-convex-combination-direction} limit the number of active piecewise intervals, and Constraints \eqref{equation:milp-convex-combination-direction1} and \eqref{equation:milp-convex-combination-direction2} allow only convex combinations of values surrounding the active interval.

\subsubsection{Linearization of Pump Head Gain Constraints}
\label{section:linearization-head-gain}%
Following the derivation of Constraints \eqref{equation:milp-convex-combination-pipe-loss}--\eqref{equation:milp-convex-combination-relations}, the piecewise \emph{lower bounds} for head gain relaxations are stated as
\begin{equation}
    g_{ij}^{k} \geq \sum_{p = 1}^{\mathclap{\lvert \mathcal{Q}_{ij}^{k +} \rvert}} (a_{ij} + b_{ij} [\hat{q}_{ijp}^{k +}]^{c_{ij}}) \lambda_{ijp}^{k +}, \, \forall (i, j) \in \mathcal{P}, \, \forall k \in \mathcal{K} \label{equation:milp-convex-combination-pump-loss}.
\end{equation}
Directed flow variables are similarly constrained by the convex combination, i.e.,
\begin{equation}
    q_{ij}^{k +} = \sum_{p = 1}^{\mathclap{\lvert \mathcal{Q}_{ij}^{k +} \rvert}} \hat{q}_{ijp}^{k +} \lambda_{ijp}^{k +}, \, \forall (i, j) \in \mathcal{P}, \, \forall k \in \mathcal{K}.
\label{equation:milp-convex-combination-pump-flow}%
\end{equation}
Activations of status-dependent convex combination variables are limited by constraints
\begin{subequations}
\begin{align}
    \sum_{p = 1}^{\mathclap{\lvert \mathcal{Q}_{ij}^{k +} \rvert}} \lambda_{ijp}^{k +} = z_{ij}, ~ \sum_{p = 2}^{\mathclap{\lvert \mathcal{Q}_{ij}^{k +} \rvert}} x_{ijp}^{k +} = z_{ij}, \, &\forall (i, j) \in \mathcal{P}, \, \forall k \in \mathcal{K} \label{equation:milp-convex-combination-activation-pump} \\
    \lambda_{ijp}^{k +} \leq x_{ijp}^{k +} + x_{ij, p + 1}^{k +}, \; &\forall (i, j) \in \mathcal{P}, \, \forall k \in \mathcal{K}, \, \forall p \in \{2, 3, \dots, \lvert \mathcal{Q}_{ij}^{k +} \rvert - 1\} \\
    \lambda_{ij, 1}^{k +} \leq x_{ij, 2}^{k +}, \; \lambda_{ij, \lvert \mathcal{Q}_{ij}^{k +} \rvert}^{k +} \leq x_{ij, \lvert \mathcal{Q}_{ij}^{k+} \rvert}^{k +}, \; &\forall (i, j) \in \mathcal{P}, \,  \forall k \in \mathcal{K} \label{equation:milp-convex-combination-pump-relations}.
\end{align}
\label{equation:pw-pump-lower-bound}%
\end{subequations}
Similar to Constraints \eqref{equation:milp-convex-combination-relations}, Constraints \eqref{equation:pw-pump-lower-bound} limit the number of active piecewise intervals and allow only convex combinations of values surrounding the active interval.
Note that convex combination variables are limited by pump statuses, $z_{ij}^{k}$, in Constraints \eqref{equation:milp-convex-combination-activation-pump}.

\subsubsection{Piecewise Relaxation of the Optimal Water Flow Problem}
Given the previously-described constraints, the MILP \emph{piecewise relaxation} of the OWF is written as
\begin{equation}
\tag{MILP-PW}
\begin{aligned}
    & \text{minimize}
    & & \textnormal{Cost of electricity} && f(q, z) ~ \textnormal{of Equation} ~ \eqref{equation:owf-objective-function} \\
    & \text{subject to}
    & & \textnormal{\eqref{equation:mincp} subset} && \textnormal{Constraints} ~ \eqref{equation:node-head-bounds}, \eqref{equation:tank-volume-integration}, \eqref{equation:mincp-valve-flow-bounds}\textnormal{--}\eqref{equation:mincp-pump-flow-bounds}, \eqref{equation:mincp-pump-head}\textnormal{--}\eqref{equation:pump-switch} \\
    & & & \textnormal{Relaxed physics} && \textnormal{Constraints} ~ \eqref{equation:micpr-direction-bounds}\textnormal{--}\eqref{equation:milp-head-loss}, \eqref{equation:milp-head-gain},  \eqref{equation:milp-convex-combination-pipe-loss}\textnormal{--}\eqref{equation:pw-pump-lower-bound}.
\end{aligned}
\label{equation:milp-pw}%
\end{equation}
where the definitions of $q_{ij}^{k}$, $(i, j) \in \mathcal{L}$, $k \in \mathcal{K}$ from Equations \eqref{equation:flow-equations} are used where necessary.
Compared to \eqref{equation:milp-oa}, \eqref{equation:milp-pw} is capable of approximating \eqref{equation:mincp} more accurately, albeit at the cost of additional binary variables that model piecewise segments.

\section{Network-based Valid Inequalities}
Section \ref{section:problem_formulation} introduced the OWF as a MINLP, and Section \ref{section:convex_relaxations} presented polyhedral relaxations that strengthen features of recent formulations encountered in the literature.
In this section, we state two sets of network-based valid inequalities (or cuts) to further strengthen these relaxations.
Section \ref{subsection:parallel-pumps} restates existing cuts for groups of identical pumps.
Section \ref{section:feasibility-step} states novel inequalities derived using duality-based physical constraint reformulations.

\label{section:valid-inequalities}
\subsection{Valid Inequalities for Parallel Identical Pumps}
\label{subsection:parallel-pumps}
In some WDNs, pumps with identical properties are installed in parallel.
The symmetry of possible activations in these pump groups can add unnecessary combinatorial complexity.
As suggested by \citet{gleixner2012towards}, the symmetry-breaking constraints are stated as 
\begin{equation}
    z_{i_{1}j_{1}}^{k} \leq z_{i_{2}j_{2}}^{k} \leq \dots \leq z_{i_{n}j_{n}}^{k}, \, \forall \{(i_{1}, j_{1}), (i_{2}, j_{2}), \dots (i_{n}, j_{n})\} \in \mathcal{P}_{G}, \, \forall k \in \mathcal{K}, \label{equation:mincp-pump-groups}
\end{equation}
where $\mathcal{P}_{G}$ comprises subsets that index identical pumps installed in parallel, which are ordered lexicographically by arbitrary indices (e.g., $1, 2, \dots, n$, above).
For each such pump group, this reduces the number of feasible decisions for active pumps from $2^{n}$ to $n + 1$.

\subsection{Novel Duality-based Valid Inequalities}
\label{section:feasibility-step}
Here, we extend a recent duality-based constraint reformulation used for optimal gravity-fed pipe network design \citep{tasseff2020exact}.
We emphasize that the inequalities we propose are for \emph{operational} WDN problems, which includes tanks, pumps, and valves.
\emph{This is a primary theoretical contribution of our paper}, which was initially proposed by the first author in their dissertation \citep{tasseff2021optimization}.
Similar inequalities have also been independently derived and explored by \citet{demassey2021optimizing}.
For brevity, we forgo derivation of these inequalities, although a detailed discussion appears in the supplemental material.
In summary, the inequalities are derived by (i) extending the convex ``Content Model'' of \citet{collins1978solving} used for pipe network analysis, (ii) relating the objective of the Content Model with the objective of its dual via strong duality, and (iii) relaxing bilinear terms of the proposed ``strong duality'' inequality involving tank heads and flows.
Letting $w_{i}^{k}$, $i \in \mathcal{T}$, $k \in \mathcal{K}$, denote McCormick-relaxed products of tank heads and flows, these inequalities are
\begin{equation}
    \begin{gathered}
        \sum_{\mathclap{(i, j) \in \mathcal{A}}} L_{ij} r_{ij} \left[(q_{ij}^{k +})^{1 + \alpha} + (q_{ij}^{k -})^{1 + \alpha}\right]
       - \sum_{\mathclap{(i, j) \in \mathcal{P}}} \left[a_{ij} q_{ij}^{k +} + b_{ij} (q_{ij}^{k +})^{c_{ij} + 1}\right]
       - \sum_{\mathclap{i \in \mathcal{T}}} w_{i}^{k} \\
       - \sum_{i \in \mathcal{R}} h_{i}^{k} \left[\sum_{(i, j) \in \delta^{k +}_{i}} (q_{ij}^{k +} - q_{ij}^{k -}) - \sum_{(j, i) \in \delta^{k -}_{i}} (q_{ji}^{k +} - q_{ji}^{k -})\right]
       - \sum_{i \in \mathcal{D}} h_{i} \overline{q}_{i}^{k} \leq 0, \, \forall k \in \mathcal{K}.
    \end{gathered}
\label{equation:convex_inequality}%
\end{equation}
Note that these new inequalities are \emph{convex} and can be applied for all time indices $k \in \mathcal{K}$.
Similar to \citet{tasseff2020exact}, these inequalities relate frictional losses, energy gains, flow production, and flow consumption to model \emph{energy conservation} at each time index $k \in \mathcal{K}$.

\section{Bound Tightening and Optimization-based Valid Inequalities}
\label{section:strengthening}
The relaxations and inequalities developed in prior sections depend strongly on the tightness of variable bounds, especially pipe and pump flow bounds.
In this section, we describe two methods: (i) improvement of variable bounds and (ii) \emph{computation} of inequalities that strengthen a relaxation.
To these ends, Section \ref{section:obbt} develops an OBBT algorithm that can leverage direction-based OWF relaxations, and Section \ref{subsection:optimization-based-valid-inequalities} extends the intuition of OBBT to develop an algorithm capable of discovering pairwise valid inequalities for the OWF.

\subsection{Optimization-based Bound Tightening}
\label{section:obbt}
\begin{algorithm}[t]
  \caption{Optimization-based bound tightening for water network problems.}
  \label{algorithm:obbt}
  \begin{algorithmic}[1]
    \Input Any direction-based model formulation, which comprises some feasible set $\Omega$
    \Output Bounds for direction-based formulations: $\underline{h}$, $\overline{h}$, $\underline{q}^{+}$, $\overline{q}^{+}$, $\underline{q}^{-}$, $\overline{q}^{-}$, $\underline{y}$, $\overline{y}$, $\underline{z}$, $\overline{z}$
    \Repeat \label{line:obbt-repeat}
    \State $(\underline{h}^{f}, \overline{h}^{f}, \underline{q}^{f + }, \overline{q}^{f + }, \underline{q}^{f -}, \overline{q}^{f -}, \underline{y}^{f}, \overline{y}^{f}, \underline{z}^{f}, \overline{z}^{f}) \gets (\underline{h}, \overline{h}, \underline{q}^{+}, \overline{q}^{+}, \underline{q}^{-}, \overline{q}^{-}, \underline{y}, \overline{y}, \underline{z}, \overline{z})$ \label{line:obbt-initialize-bounds}
        \State $\Omega \gets$ Direction-based formulation given $(\underline{h}^{f}, \overline{h}^{f}, \underline{q}^{f + }, \overline{q}^{f + }, \underline{q}^{f -}, \overline{q}^{f -}, \underline{y}^{f}, \overline{y}^{f}, \underline{z}^{f}, \overline{z}^{f})$ \label{line:obbt-make-relaxation}
        \ForAll{$(i \in \mathcal{D} \cup \mathcal{R}, \, k \in \mathcal{K})$, $(i \in \mathcal{T}, \, k \in \tilde{\mathcal{K}})$} \label{line:obbt-head-bounds-1}
            \State $(\underline{h}_{i}^{k}, \overline{h}_{i}^{k}) \gets (\textrm{minimize} ~ h_{i}^{k} ~ \textrm{s.t.} ~ \Omega, ~ \textrm{maximize} ~ h_{i}^{k} ~ \textrm{s.t.} ~ \Omega)$
        \EndFor \label{line:obbt-head-bounds-4}
        \ForAll{$(i, j) \in \mathcal{A} \cup \mathcal{V}, \, k \in \mathcal{K}$} \label{line:obbt-flow-bounds-1}
        \State $(\underline{q}_{ij}^{k \pm}, \overline{q}_{ij}^{k \pm}) \gets (\textrm{minimize} ~ q_{ij}^{k \pm} ~ \textrm{s.t.} ~ \Omega \cup \{y_{ij}^{k \pm} = 1\}, ~ \textrm{maximize} ~ q_{ij}^{k \pm} ~ \textrm{s.t.} ~ \Omega \cup \{y_{ij}^{k \pm} = 1\})$ \label{line:obbt-flow-bounds-2}
        \State $(\underline{y}_{ij}^{k}, \overline{y}_{ij}^{k}) \gets (\textrm{minimize} ~ y_{ij}^{k} ~ \textrm{s.t.} ~ \Omega, ~ \textrm{maximize} ~ y_{ij}^{k} ~ \textrm{s.t.} ~ \Omega)$ \label{line:obbt-flow-bounds-2-1}
        \EndFor \label{line:obbt-flow-bounds-4}
        \ForAll{$(i, j) \in \mathcal{P} \cup \mathcal{V}, \, k \in \mathcal{K}$} \label{line:obbt-status-bounds-1}
        \State $(\underline{z}_{ij}^{k}, \overline{z}_{ij}^{k}) \gets (\textrm{minimize} ~ z_{ij}^{k} ~ \textrm{s.t.} ~ \Omega, ~ \textrm{maximize} ~ z_{ij}^{k} ~ \textrm{s.t.} ~ \Omega)$ \label{line:obbt-status-bounds-2}
        \EndFor \label{line:obbt-status-bounds-4}
        \ForAll{$(i, j) \in \mathcal{P}, \, k \in \mathcal{K}$} \label{line:obbt-epsilon-bounds-1}
        \State $(\underline{q}_{ij}^{k +}, \overline{q}_{ij}^{k +}) \gets (\textrm{minimize} ~ q_{ij}^{k +} ~ \textrm{s.t.} ~ \Omega \cup \{z_{ij}^{k} = 1\}, ~ \textrm{maximize} ~ q_{ij}^{k +} ~ \textrm{s.t.} ~ \Omega \cup \{z_{ij}^{k} = 1\})$ \label{line:obbt-epsilon-bounds-2}
        \EndFor \label{line:obbt-epsilon-bounds-3}
    \Until{$(\underline{h}^{f}, \overline{h}^{f}, \underline{q}^{f + }, \overline{q}^{f + }, \underline{q}^{f -}, \overline{q}^{f -}, \underline{y}^{f}, \overline{y}^{f}, \underline{z}^{f}, \overline{z}^{f}) = (\underline{h}, \overline{h}, \underline{q}^{+}, \overline{q}^{+}, \underline{q}^{-}, \overline{q}^{-}, \underline{y}, \overline{y}, \underline{z}, \overline{z})$} \label{line:obbt-repeat-end}
  \end{algorithmic}
\end{algorithm}

Algorithm \ref{algorithm:obbt} presents the variant of the OBBT algorithm used in this paper.
Here, the input $\Omega$ comprises the variable bounds and constraints used in some direction-based formulation of the OWF or a relaxation thereof, e.g., $\Omega = \eqref{equation:milp-oa}$ or $\Omega = \eqref{equation:milp-pw}$.
Lines \ref{line:obbt-repeat} and \ref{line:obbt-repeat-end} give the algorithm's iteration and termination conditions, respectively, i.e., repeat the interior of the algorithm (Lines \ref{line:obbt-initialize-bounds}--\ref{line:obbt-epsilon-bounds-3}) until variable bounds no longer improve.
Note that this is an idealized termination condition, and actual termination depends on the desired bound precision, maximum number of iterations, or convergence tolerance.
Line \ref{line:obbt-initialize-bounds} sets the variable bounds for the current iteration of the algorithm.
Line \ref{line:obbt-make-relaxation} constructs a new relaxation, $\Omega$, using these variable bounds.
Lines \ref{line:obbt-head-bounds-1}--\ref{line:obbt-head-bounds-4} derive new bounds for heads in the network.
Here, Line \ref{line:obbt-head-bounds-1} solves a minimization problem that yields a lower bound, $\underline{h}_{i}^{k}$, as well as a maximization problem that yields an upper bound, $\overline{h}_{i}^{k}$.
In a similar manner, Lines \ref{line:obbt-flow-bounds-1}--\ref{line:obbt-flow-bounds-4} compute new lower and upper bounds for pipe and valve flow and direction variables in the network, and Lines \ref{line:obbt-status-bounds-1}--\ref{line:obbt-status-bounds-4} compute new bounds for pump and valve status variables.
Finally, Lines \ref{line:obbt-epsilon-bounds-1}--\ref{line:obbt-epsilon-bounds-3} compute new bounds for variable pump flows \emph{when the pump is active} (i.e., $z_{ij}^{k} = 1$).
This conditional bound tightening notably allows for further strengthening of relaxations for each pump's head gain curve by improving the fixed $\underline{q}_{ij}^{k +}$.
We finally emphasize that, as later elaborated upon in Section \ref{subsection:effects_of_obbt}, the variations of this OBBT algorithm we ultimately employ frequently solve subproblems that impose integrality constraints.

\subsection{Optimization-based Valid Inequalities}
\label{subsection:optimization-based-valid-inequalities}
\begin{algorithm}[t]
\caption{Binary-binary OBCG procedure, where one variable is fixed to zero.}
\label{algorithm:obcg-1}
\begin{algorithmic}[1]
  \Input Any direction-based model formulation, which comprises some feasible set $\Omega$
  \Output Valid inequalities for direction-based model formulations, denoted by $\bar{\mathcal{X}}$
  \State
  $\mathcal{B} \gets \{y_{ij}^{k} : (i, j) \in \mathcal{L}\} \cup \{z_{ij}^{k} : (i, j) \in \mathcal{P} \cup \mathcal{V}\}$, $\bar{\mathcal{X}} = \emptyset$ \label{line:obcg-1-collect-variable-sets}
  \ForAll{$(x_{1}, x_{2}) \in (\mathcal{B} \times \mathcal{B}) \setminus \{(x, x) : x \in \mathcal{B}\}$} \label{line:obcg-1-binary-binary-loop}
      \State $(\underline{x}_{1}^{0}, \overline{x}_{1}^{0}) \gets (\textrm{minimize} ~ x_{1} ~ \textrm{s.t.} ~ \Omega \cup \{x_{2} = 0\}, ~ \textrm{maximize} ~ x_{1} ~ \textrm{s.t.} ~ \Omega \cup \{x_{2} = 0\})$ \label{line:obcg-1-binary-binary-minimize-0}
      \If{$\underline{x}_{1}^{0} = \overline{x}_{1}^{0} = 0$} \label{line:obcg-1-if-binary-binary-0-0}
          \State $\bar{\mathcal{X}} \gets \bar{\mathcal{X}} \cup \{x_{1} \leq x_{2}\}$ \label{line:obcg-1-cut-binary-binary-0-0}
      \ElsIf{$\underline{x}_{1}^{0} = \overline{x}_{1}^{0} = 1$} \label{line:obcg-1-if-binary-binary-0-1}
          \State $\bar{\mathcal{X}} \gets \bar{\mathcal{X}} \cup \{x_{1} + x_{2} \geq 1\}$ \label{line:obcg-1-cut-binary-binary-0-1}
      \EndIf \label{line:obcg-1-endif-binary-binary-0}
  \EndFor
\end{algorithmic}
\end{algorithm}

Algorithm \ref{algorithm:obbt} shows how some variable bounds (e.g., the minimum flow rate through an \emph{active} pump, $\underline{q}_{ij}^{k +}$) can be improved for \emph{continuous} variables (e.g., $q_{ij}^{k +}$) that are \emph{conditional} upon the value of \emph{discrete} variables (e.g., when $z_{ij}^{k} = 1$).
In this subsection, we broaden this intuition to compute \emph{conditional} valid inequalities between arbitrary pairs of variables.
Algorithm \ref{algorithm:obcg-1} outlines one of these novel OBCG algorithms.
Here, Line \ref{line:obcg-1-collect-variable-sets} instantiates the set of binary variables to consider in the cut generation procedure.
Line \ref{line:obcg-1-binary-binary-loop} defines the binary-binary cut generation loop, where all unique, ordered pairs of binary variables are considered.
Line \ref{line:obcg-1-binary-binary-minimize-0} minimizes and maximizes the first variable in the pair, respectively, subject to the relaxation constraints and a fixing of the second variable to zero.
On Line \ref{line:obcg-1-if-binary-binary-0-0}, if both minimization and maximization imply objectives of zero, a cut relating the two discrete variables can be derived, which is added to the set of cuts $\bar{\mathcal{X}}$ on Line \ref{line:obcg-1-cut-binary-binary-0-0}.
On the other hand, if minimization and maximization both yield one, a similar process is followed on Lines \ref{line:obcg-1-if-binary-binary-0-1}--\ref{line:obcg-1-endif-binary-binary-0}.
Otherwise, if the minimization and maximization problems yield different values, we do not add any cuts relating this pair.
In another algorithm included in the supplemental material, the same process is repeated while fixing the second variable to one instead of zero.
Both algorithms generate valid binary-binary inequalities for $k \in \mathcal{K}$.

\begin{algorithm}[t]
\caption{Binary-continuous OBCG procedure, where the binary variable is fixed.}
\label{algorithm:obcg-3}
\begin{algorithmic}[1]
  \Input Any direction-based model formulation, which comprises some feasible set $\Omega$
  \Output Valid inequalities for direction-based model formulations, denoted by $\bar{\mathcal{X}}$
  \ForAll{$(x_{1}, x_{2}) \in \left(\{q_{ij}^{k} : (i, j) \in \mathcal{L}\} \cup \{q_{i}^{k} : i \in \mathcal{T}\} \cup \{h_{i}^{k} : i \in \mathcal{N}\}\right) \times \mathcal{B}$} \label{line:obcg-3-continuous-binary-loop}
      \State $\underline{x}_{1}^{0} \gets \textrm{minimize} ~ x_{1} ~ \textrm{s.t.} ~ \Omega \cup \{x_{2} = 0\}$ \label{line:obcg-3-continuous-binary-minimize-0}
      \State $\bar{\mathcal{X}} \gets \bar{\mathcal{X}} \cup \{\underline{x}_{1}^{0} (1 - x_{2}) + \underline{x}_{1} x_{2} \leq x_{1}\}$ \label{line:obcg-3-continuous-binary-minimize-cut-0}
      \State $\overline{x}_{1}^{0} \gets \textrm{maximize} ~ x_{1} ~ \textrm{s.t.} ~ \Omega \cup \{x_{2} = 0\}$ \label{line:obcg-3-continuous-binary-maximize-0}
      \State $\bar{\mathcal{X}} \gets \bar{\mathcal{X}} \cup \{x_{1} \leq \overline{x}_{1}^{0} (1 - x_{2}) + \overline{x}_{1} x_{2}\}$
      \State $\underline{x}_{1}^{1} \gets \textrm{minimize} ~ x_{1} ~ \textrm{s.t.} ~ \Omega \cup \{x_{2} = 1\}$
      \State $\bar{\mathcal{X}} \gets \bar{\mathcal{X}} \cup \{\underline{x}_{1}^{1} x_{2} + \underline{x}_{1} (1 - x_{2}) \leq x_{1}\}$
      \State $\overline{x}_{1}^{1} \gets \textrm{maximize} ~ x_{1} ~ \textrm{s.t.} ~ \Omega \cup \{x_{2} = 1\}$
      \State $\bar{\mathcal{X}} \gets \bar{\mathcal{X}} \cup \{x_{1} \leq \overline{x}_{1}^{1} x_{2} + \overline{x}_{1} (1 - x_{2})\}$ \label{line:obcg-3-continuous-binary-maximize-cut-1}
  \EndFor
\end{algorithmic}
\end{algorithm}

Algorithm \ref{algorithm:obcg-3} outlines a variant of OBCG where variable pairs comprise \emph{continuous} variables of flow and head as well as \emph{binary} variables of flow direction and control component status.
The goal of this algorithm is to derive cuts that improve continuous variable bounds \emph{depending on the value of each discrete variable}.
For example, Line \ref{line:obcg-3-continuous-binary-minimize-0} minimizes a continuous variable subject to the relaxation constraints and a fixing of the binary variable to zero.
In turn, the optimal solution provides a potentially tighter \emph{lower} bound for the continuous variable \emph{when the binary variable is zero}.
This inequality is added to the set of cuts $\bar{\mathcal{X}}$ on Line \ref{line:obcg-3-continuous-binary-minimize-cut-0}.
Similar cuts are derived for other bounds and fixings on Lines \ref{line:obcg-3-continuous-binary-maximize-0}--\ref{line:obcg-3-continuous-binary-maximize-cut-1}.

We remark that the OBCG algorithms can be computationally expensive due to quadratic scaling in the number of variables.
In this paper, we thus restrict each algorithm's execution to an individual steady state $k \in \mathcal{K}$, where Constraints \eqref{equation:tank-volume-integration}, \eqref{equation:tank-volume-recovery}, and \eqref{equation:pump-switch} are omitted.
This implies that \emph{our OBCG methods compute inequalities that are valid for each steady state and thus also valid for the original, intertemporal OWF problem}.
In Section \ref{section:computational_experiments}, we explore similar relaxations of the OWF for the OBBT algorithm of Section \ref{section:obbt}.

\section{Computational Experiments}
\label{section:computational_experiments}
This section empirically evaluates the contributions of our paper.
Section \ref{subsection:experimental_setup} describes the software, instances, resources, and parameters used in the experiments.
Section \ref{subsection:effects_of_obbt} explores the impacts of OBBT on the root relaxation of \eqref{equation:milp-oa}, showing that improved OBBT methods lead to substantially improved dual bounds.
Section \ref{subsection:effects_of_almost_exact_inequalities} explores the same for the duality-based inequalities of Section \ref{section:feasibility-step}, showing significant relaxation improvements on challenging OWF instances.
Section \ref{subsection:effects_of_valid_inequalities} evaluates the OBCG inequalities using the algorithms of Section \ref{subsection:optimization-based-valid-inequalities}, which result in additional improvements.
Section \ref{subsection:effects_of_polyhedral_parameterizations} evaluates the primal bound quality from various polyhedral relaxations, showing an interesting tradeoff between combinatorial complexity and the production of \eqref{equation:mincp}-feasible solutions.
Finally, Section \ref{subsection:relative-effects} evaluates the relative effects of all algorithmic contributions on primal and dual bounds for a set of large, complex OWF instances, showing that \emph{the combination} of our contributions is necessary to achieve increased qualities of primal and dual bounds.

\subsection{Experimental Setup}
\label{subsection:experimental_setup}
\begin{table}[t]
    \small
    \begin{center}
        \begin{tabular}{|c|cccccccccccc|}
            \hline
            Network & $\lvert \mathcal{N} \rvert$ & $\lvert \mathcal{D} \rvert$ & $\lvert \mathcal{R} \rvert$& $\lvert \mathcal{T} \rvert$ & $\lvert \mathcal{L} \rvert$ & $\lvert \mathcal{A} \rvert$ & $\lvert \mathcal{P} \rvert$ & $\lvert \mathcal{V} \rvert$ & $K$ & $N_{ij}$ & $\tau^{\textrm{on}}$ & $\tau^{\textrm{off}}$ \\ \hline
            \texttt{Simple FSD} & $6$  & $2$  & $3$ & $1$ & $5$  & $2$  & $3$ & $0$ & $\{12, 24, 48\}$ & $18$ & $3{,}600$ s & $1{,}800$ s \\
            \texttt{AT(M)}      & $24$ & $19$ & $3$ & $2$ & $44$ & $41$ & $3$ & $0$ & $\{12, 24, 48\}$ & $18$ & $3{,}600$ s & $1{,}800$ s \\
            \texttt{Poormond}   & $52$ & $46$ & $1$ & $5$ & $55$ & $44$ & $7$ & $4$ & $\{12, 24, 48\}$ & $6$ & $3{,}600$ s & $1{,}800$ s \\ \hline
        \end{tabular}
    \end{center}
    \caption{Summary of the OWF instances considered in this paper using the notation from Section \ref{section:network}, e.g., where $\mathcal{N}$ is the set of nodes, $\mathcal{L}$ is the set of node-connecting components, and $K$ is the number of time steps.}
    \label{table:network_summary}
\end{table}

All formulations and algorithms were implemented in the \textsc{Julia} programming language using \textsc{JuMP}, version 0.22 \citep{dunning2017jump}, and version 0.9 of \textsc{WaterModels}, an open-source \textsc{Julia} package for WDN optimization \citep{watermodels}.
Hereafter, a parameter, $\xi$, is used to concisely describe the flow partitioning of a MILP relaxation, where \emph{smaller} values of $\xi$ indicate a \emph{greater} number of outer- and/or inner-approximation points.
Specifically, $\xi$ represents the maximum error tolerance (in meters) of a flow partitioning scheme when employing \eqref{equation:milp-pw} for modeling pipe head loss and pump head gain nonlinearities.
That is, $\xi$ is a threshold for the distance between purple and orange lines in Figure \ref{figure:polyhedral_relaxations}.
In practice, this parameter is used to construct the partitions $\mathcal{Q}_{ij}^{\pm}$ for all pipes and pumps via the \textsc{PolyhedralRelaxations} \textsc{Julia} package of \citet{sundar_sequence_2021}.

The experiments consider three networks of various sizes that appear in the pump scheduling literature and are summarized in Table \ref{table:network_summary} \citep{bonvin2021pump,mala2018lost}.
Visualizations of the networks are presented in the supplemental material.
Although \texttt{Simple FSD} and \texttt{Poormond} are similar (but not equivalent) to the networks studied by \cite{bonvin2021pump}, our version of \texttt{AT(M)} increases maximum tank water levels from $71.53$ to $72.00$ meters.
Finally, all pumps in \texttt{Simple FSD} and \texttt{AT(M)} are identical and leverage the symmetry-breaking Constraints \eqref{equation:mincp-pump-groups}, but all pumps in \texttt{Poormond} are unique.

Each of the three networks have correspondence with fifteen unique OWF instances, which were also derived from \cite{bonvin2021pump}.
Each instance varies in two ways.
First, although each instance assumes a twenty-four hour scheduling period, the patterns of demands and electricity prices are specified in two hour, one hour, and thirty minute intervals.
That is, instances have various temporal dimensions, where $K \in \{12, 24, 48\}$.
Second, for each $K$, electricity price patterns differ across five daily estimates.
Temporal profiles of these price patterns are illustrated in the supplemental material.
In summary, the instances study variations in temporal resolution of the model and in electricity pricing.

Each experiment was executed on a node containing two AMD EPYC 7H12 processors, providing a total of 128 cores at 2.60 GHz, and 512 GB of memory.
For the solution of all MILPs, \textsc{Gurobi} 9.1 was used with the parameter \texttt{MIPGap=0.0}.
In the execution of OBBT and OBCG algorithms, subproblems were provided a wall-clock limit of one minute.
\emph{If exceeded, the final dual bound of the subproblem was used in place of its optimal solution.}

\subsection{Effects of Optimization-based Bound Tightening}
\label{subsection:effects_of_obbt}
\begin{table}[t]
    \small
    \begin{center}
        \begin{tabular}{|c|cccp{7cm}|}
            \cline{1-5}
            Name & Formulation & Cont. Rel. & Time Cmplx. & Notes \\
            \cline{1-5}
            BT-SR & \eqref{equation:milp-oa} & \ding{51} & Reduced  & Relaxes $\mathcal{D}$, $\mathcal{R}$; ignores intertemporal constrs. \\
            BT-SS & \eqref{equation:milp-pw} & \ding{55} & Reduced  & Relaxes $\mathcal{D}$, $\mathcal{R}$; ignores intertemporal constrs. \\
            BT-SQ & \eqref{equation:milp-pw} & \ding{55} & Sequence & Ignores intertemporal constraints \\
            BT-TR & \eqref{equation:milp-oa} & \ding{51} & Complete & Continuous relaxation of \eqref{equation:milp-oa} \\
            BT-TS & \eqref{equation:milp-pw} & Partial   & Complete & Continuous relaxation of $k \neq k^{\prime} \in \mathcal{K}$ \\
            \cline{1-5}
        \end{tabular}
        \caption{Summary of the various optimization-based bound tightening methods considered in this paper.}
        \label{table:obbt_summary}%
    \end{center}
\end{table}

This section evaluates various OBBT algorithms inspired by Section \ref{section:obbt} and their effects on the root relaxation of \eqref{equation:milp-oa}.
Specifically, \emph{we show that different OBBT techniques, which vary in computational complexity, can provide significant relaxation improvements, albeit at different computational costs.}
Each OBBT algorithm was provided a wall-clock time limit of eight hours or a maximum of $25$ iterations.
The subproblems for each algorithm were trivially parallelized over $128$ threads available on one processor of a compute node.
Each relaxation employed within an algorithm assumed an error tolerance of $\xi = 1$ meter.

We consider five different variants of OBBT, each with greater temporal complexity than the last, and each of which begins with bounds derived from the last.
These methods and their naming conventions are summarized in Table \ref{table:obbt_summary}.
Here, BT-SR assumes a \emph{single-step} approximation of \eqref{equation:mincp} that ignores the intertemporal Constraints \eqref{equation:tank-volume-integration} and where variables are bounded as $\underline{x}_{i} = \min_{k \in \mathcal{K}} \{\underline{x}_{i}^{k}\}$ and $\overline{x}_{i} = \max_{k \in \mathcal{K}} \{\overline{x}_{i}^{k}\}$.
Since this is a relaxation of the constraints for \emph{any steady state} $k \in \mathcal{K}$, the bounds determined via BT-SR can be applied for \emph{any} $k \in \mathcal{K}$.
This method continuously relaxes a \eqref{equation:milp-oa}-type relaxation.

Similar to BT-SR, BT-SS uses the above assumptions but instead leverages a \eqref{equation:milp-pw}-type relaxation \emph{without relaxing integrality}.
Next, BT-SQ performs OBBT for the \emph{sequence of steady states} $k \in \mathcal{K}$, i.e., there is no temporal relaxation of variable bounds, although intertemporal constraints are ignored.
Finally, BT-TR and BT-TS perform OBBT for the complete OWF problem.
Here, BT-TR assumes a continuous relaxation of \eqref{equation:milp-oa}, while BT-TS is a partial continuous relaxation of \eqref{equation:milp-pw}.
Specifically, if a variable at step $k^{\prime} \in \mathcal{K}$ is tightened, integrality is enforced for $k^{\prime}$ but relaxed for $k \neq k^{\prime}$.

\input{figure_obbt_root_relaxation.tex}

Figure \ref{figure:obbt-improvement} illustrates improvements in each root relaxation's objective for \eqref{equation:milp-oa}, assuming $\xi = 1$ meter, and after applying each of the five bound tightening techniques.
For these experiments, \textsc{Gurobi} 9.1 was used with the additional parameters \texttt{NumericFocus=3}, \texttt{FeasibilityTol=1.0e-9}, and \texttt{OptimalityTol=1.0e-9} to ensure highly accurate estimates of each root relaxation's objective.
Here, each improvement is measured in comparison with a $\xi = 1$ meter \eqref{equation:milp-oa} relaxation that assumes naive, network-based bounds, i.e.,
\begin{equation}
    \textrm{\% Improvement} = 100\% \left(\frac{f_{2} - f_{1}}{f_{1}}\right),
\end{equation}
where $f_{1}$ and $f_{2}$ are objectives when using naive and OBBT-based bounds, respectively.

To begin, Figure \ref{figure:obbt-improvement} indicates that the effects of OBBT are often dependent on the network topology.
Whereas OBBT generally provides improvements greater than $25\%$ for \texttt{Simple FSD} and \texttt{Poormond} instances, smaller improvements are observed for \texttt{AT(M)}.
This could be because flow directions are generally predetermined for the former two networks but not the latter, which has a large number of loops and relatively large pipe flow bounds.
The second observation is that BT-SR, the weakest and most efficient technique, often provides most of the improvement over naive bounds.
On the other hand, BT-SS provides generally small improvements over BT-SR, although the improvements are larger for the more complex instances of \texttt{AT(M)} and \texttt{Poormond}.
Similarly, BT-SQ provides small but potentially important improvements over BT-SS.
Interestingly, BT-TR provides mixed improvements over all instances, but sometimes gives rise to extreme improvements, as exemplified in the last \texttt{Poormond} instance for $K = 12$.
Finally, BT-TS is observed to provide almost no improvement over BT-TR, suggesting BT-TS is likely not worth its computational cost.

\input{figure_obbt_timing.tex}

Figure \ref{figure:obbt-timing} plots the maximum wall clock time of each OBBT technique for the various networks and time discretizations.
First, the execution times for the algorithms often depend on the network size and time discretization.
For \texttt{Simple FSD}, all OBBT routines complete within a few minutes, although techniques that model the full time series (BT-TR and BT-TS) require considerably more time than others.
For \texttt{AT(M)}, some techniques require hours to complete, and full time series techniques could be viewed as prohibitively expensive.
For \texttt{Poormond}, the same trends are observed, although BT-SQ completes within minutes.
We remark, however, that the OBBT methods are amenable to parallelization, and given sufficient computational resources, an idealized parallel version of the most expensive algorithm would complete within nine minutes.
Nonetheless, Figures \ref{figure:obbt-improvement} and \ref{figure:obbt-timing} emphasize there are tradeoffs between the improvement of variable bounds and the required resources.

\subsection{Effects of Duality-based Inequalities}
\label{subsection:effects_of_almost_exact_inequalities}
This section evaluates the efficacy of the duality-based inequalities of Section \ref{section:feasibility-step}.
Specifically, \emph{we show that duality-based inequalities can significantly strengthen OWF relaxations}.
Further, these inequalities can be included within an OWF formulation without the need for preprocessing, as compared to the previous OBBT techniques.
Here, the inequalities were implemented as linear relaxations of Constraints \eqref{equation:convex_inequality}, where nonlinear terms involving $q^{\pm}$ were instead linearly outer-approximated using the partition sets $\mathcal{Q}^{\pm}$, as with the head loss and head gain nonlinearities of \eqref{equation:milp-oa}.
We forgo an explicit derivation and statement of the linear relaxations of these constraints for brevity, although they bear great similarity with the linearizations described by \cite{tasseff2020exact} and Section \ref{section:mixed-integer-linear-relaxation}.

\input{figure_cut_root_relaxation.tex}

Similar to Figure \ref{figure:obbt-improvement}, Figure \ref{figure:cut-improvement} shows the improvement in the root relaxation of \eqref{equation:milp-oa} with $\xi = 1$ meter after applying the various valid inequalities developed in this paper.
All experiments assumed variable bounds discovered via the final BT-TS OBBT method discussed in Section \ref{subsection:effects_of_obbt} and the same \textsc{Gurobi} parameterization as in that section's root relaxation experiments.
Each percentage improvement is measured with respect to the root relaxation objective when using BT-TS-based bounds alone, as used in the final experiments of Section \ref{subsection:effects_of_obbt}.
We finally note that Figure \ref{figure:cut-improvement} also includes improvements from the OBCG-derived cuts, although we reserve discussion of these results to the next subsection.

First, it's apparent that for \texttt{Simple FSD} and \texttt{AT(M)}, the improvements from incorporating duality-based cuts are marginal and always less than $0.2\%$.
This could be due to (i) the small size of \texttt{Simple FSD}, whose instance relaxations may already be relatively tight from incorporating the OBBT-based bounds of Section \ref{subsection:effects_of_obbt}; (ii) the relatively small numbers of tanks and unique pumps in these two networks, which play important roles in Constraints \eqref{equation:convex_inequality}; and (iii) the symmetry of pumps in these two networks, which already limit the operational possibilities substantially via Constraints \eqref{equation:mincp-pump-groups}.
Finally, for \texttt{Poormond}, the benefits are more dramatic: incorporating duality-based cuts provides improvements between $5\%$ and $11\%$, which could be substantial when attempting to prove optimality of a solution.

\subsection{Effects of Optimality-based Cut Generation Inequalities}
\label{subsection:effects_of_valid_inequalities}
Figure \ref{figure:cut-improvement} illustrates root relaxation improvements when including the inequalities generated via the OBCG methods of Section \ref{subsection:optimization-based-valid-inequalities}, on top of the inequalities explored in Section \ref{subsection:effects_of_almost_exact_inequalities}.
In the generation of these inequalities, all OBCG algorithms (e.g., Algorithm \ref{algorithm:obcg-1}) employed a \eqref{equation:milp-pw}-type relaxation using a parameterization of $\xi = 1$ meter, and all subproblems were parallelized over $128$ threads, as in the OBBT subproblems discussed in Section \ref{subsection:effects_of_obbt}.
We remark that, in practice, the OBCG algorithms employed here assumed subproblems similar to the BT-SQ method, i.e., where the intertemporal constraints of \eqref{equation:milp-pw} were ignored.
Thus, OBCG was used to compute inequalities for \emph{each steady state} $k \in \mathcal{K}$.

As observed in Figure \ref{figure:cut-improvement}, the OBCG-based cuts provide some improvements over duality-based cuts alone.
For \texttt{Simple FSD} and \texttt{AT(M)}, the improvements are small and always less than $0.5\%$.
However, for \texttt{Poormond}, these improvements are significant -- up to 10\%.
The effects are more dramatic for \texttt{Poormond} because of the greater complexity of the network and larger number of operational possibilities.
Presumably, the OBCG-based cuts remove a large number of operational possibilities that would be infeasible in practice.
Thus, OBCG-based cuts are capable of improving the root relaxation by a substantial amount.

\begin{table}[t]
    \small
    \begin{center}
        \begin{tabular}{c|cc|cc|cc|}
            \cline{2-7} & \multicolumn{2}{c|}{$K = 12$} &
                          \multicolumn{2}{c|}{$K = 24$} &
                          \multicolumn{2}{c|}{$K = 48$} \\
            \hline
            \multicolumn{1}{|c|}{Network} & Wall Time & Ideal Time & Wall Time & Ideal Time & Wall Time & Ideal Time \\ \hline
            \multicolumn{1}{|c|}{\texttt{Simple FSD}}
                & $0.4$ s        & $0.4$ s   & $6.8$ s        & $6.7$ s   & $6.0$ s        & $5.8$ s \\
            \multicolumn{1}{|c|}{\texttt{AT(M)}}
                & $4{,}031.3$ s  & $150.2$ s & $10{,}479.1$ s & $367.9$ s & $21{,}406.6$ s & $749.2$ s \\
            \multicolumn{1}{|c|}{\texttt{Poormond}}
                & $120.0$ s      & $15.6$ s  & $238.0$ s      & $32.8$ s  & $532.5$ s      & $92.3$ s \\
            \hline
        \end{tabular}
    \end{center}
    \caption{Wall and idealized parallel timing statistics for the optimization-based cut generation technique.}
    \label{table:obcg-timing}
\end{table}

Table \ref{table:obcg-timing} evaluates the resource tradeoffs of OBCG in comparison with the benefits illustrated in Figure \ref{figure:cut-improvement}.
First, for \texttt{Simple FSD} and \texttt{Poormond}, all pairwise cuts are computed within ten minutes of wall clock time, which is relatively small given the large number of OBCG subproblems.
On the other hand, \texttt{AT(M)} requires considerable time to compute these cuts, taking around six hours in the worst case.
This is likely due to the complexity of \eqref{equation:milp-pw}-type subproblems for \texttt{AT(M)}, which require a large number of binary variables to model piecewise pipe head loss constraints.
This time could be reduced by employing a weaker subproblem relaxation or larger error tolerance (e.g., $\xi = 10$ meters).
We also remark that the worst-case ideal parallel time for OBCG is around twelve minutes.
Nonetheless, compared with the sometimes small benefits from OBCG displayed in Figure \ref{figure:cut-improvement}, it could be concluded that there is an important tradeoff between the time required to generate these cuts and the benefits of using them, which may be topology- or bound-dependent.

\subsection{Primal-bounding Experiments}
\label{subsection:effects_of_polyhedral_parameterizations}
In this section, we evaluate six relaxations and measure their efficacy in determining primal bounds of the \emph{nonconvex} OWF.
The experiments also incorporate the model strengthening techniques explored in previous sections (i.e., OBBT and all valid inequalities).
The first relaxations are \eqref{equation:milp-oa} formulations where $\xi \in \{1, 5, 25\}$ meters and all direction variables, $y_{ij}^{k}$, $(i, j) \in \mathcal{L}$, $k \in \mathcal{K}$, are continuously relaxed.
The second relaxations are \eqref{equation:milp-pw} formulations where $\xi \in \{1, 5, 25\}$ meters, and where the duality-based Constraints \eqref{equation:convex_inequality} are inner-/outer-approximated via the sharing of discrete and continuous convex combination variables introduced in Sections \ref{section:linearization-head-loss} and \ref{section:linearization-head-gain}.
In summary, relaxations in the second set are tighter but more expensive than the first.
All experiments used a time limit of one hour and the \textsc{Gurobi} parameters \texttt{MIPGap=0.0}, \texttt{PreCrush=1}, \texttt{Presolve=2}, and \texttt{Threads=64}.
Additionally, progressively larger \texttt{BranchPriority} values were assigned to discrete variables that appear \emph{earlier} in the time series model, and larger priorities were assigned to component status variables ($z_{ij}$) when compared with direction variables ($y_{ij}$).

To ensure primal solutions are feasible with respect to \eqref{equation:mincp}, a feasibility-checking procedure, similar to those used by \cite{raghunathan2013global} and \cite{tasseff2020exact} for WDN design and \cite{bonvin2021pump} for WDN operation, was used within the MILP search.
Specifically, at every integer-feasible node of the search tree, the network analysis algorithm of \citet{todini1988gradient} was used for the corresponding schedule of controllable components.
As in \cite{bonvin2021pump} and \cite{naoum2015simulation}, when a bound infeasibility is discovered, the following combinatorial ``no good'' cut is then appended:
\begin{equation}
    \mathlarger{\sum}_{k = 1}^{K^{\textrm{inf}}} \left(\sum_{(i, j) \in \mathcal{P} \cup \mathcal{V} : \hat{z}_{ij}^{k} = 0} z_{ij}^{k} - \sum_{(i, j) \in \mathcal{P} \cup \mathcal{V} : \hat{z}_{ij}^{k} = 1} (1 - z_{ij}^{k})\right) \geq 1
    \label{equation:no-good-cut}.
\end{equation}
Here, $K^{\textrm{inf}} \in \mathcal{K}$ is the first time index at which infeasibility is detected in the extended period analysis, and $\hat{z}$ corresponds to the solution at the current BB node.
This cut implies that at least one component status must change to address the infeasibility.
If the integer solution is instead feasible, the solution to the extended period analysis is used to compute a \emph{valid} primal bound to \eqref{equation:mincp}, which is not necessarily the \emph{optimal} objective value.

\begin{figure}[t]
    \centering
    \begin{tikzpicture}[baseline=(current bounding box.north)]
        \begin{axis}[ylabel={$\%$ from Best},enlargelimits=false,height=4.3cm,
            legend style={at={(0.05, 0.94)},anchor=north west,font=\footnotesize},
            xlabel=Price Profile (Day),xticklabels={1,2,3,4,5},xtick={1,2,3,4,5},
            width=0.37\linewidth,xtick pos=left,ytick pos=left,xmin=1,xmax=5,
            ymin=0.0,ymax=2.25,legend columns=1,legend entries={$\xi = 1$ m, $\xi = 5$ m,
            $\xi = 25$ m,\eqref{equation:milp-oa},\eqref{equation:milp-pw}},legend cell align={left}]
            \node[anchor=north east] at (rel axis cs:0.99,0.98) {$K = 12$};
            \addlegendimage{index of colormap=0 of Dark2-8,mark=square*,only marks,opacity=0.5};
            \addlegendimage{index of colormap=1 of Dark2-8,mark=square*,only marks,opacity=0.5};
            \addlegendimage{index of colormap=2 of Dark2-8,mark=square*,only marks,opacity=0.5};
            \addlegendimage{black,mark=*,only marks};
            \addlegendimage{black,mark=diamond*,only marks};
            \pgfplotstableread[col sep = comma]{MATM-12_Steps-UB.csv}{\data};
            \addplot[index of colormap=0 of Dark2-8,mark=*,only marks,mark size=3pt,opacity=0.5] table [x=day,y=ub_lrdx_1] {\data};
            \addplot[index of colormap=1 of Dark2-8,mark=*,only marks,mark size=3pt,opacity=0.5] table [x=day,y=ub_lrdx_2] {\data};
            \addplot[index of colormap=2 of Dark2-8,mark=*,only marks,mark size=3pt,opacity=0.5] table [x=day,y=ub_lrdx_3] {\data};
            \addplot[index of colormap=0 of Dark2-8,mark=diamond*,only marks,mark size=3pt,opacity=0.5] table [x=day,y=ub_pwlrdx_1] {\data};
            \addplot[index of colormap=1 of Dark2-8,mark=diamond*,only marks,mark size=3pt,opacity=0.5] table [x=day,y=ub_pwlrdx_2] {\data};
            \addplot[index of colormap=2 of Dark2-8,mark=diamond*,only marks,mark size=3pt,opacity=0.5] table [x=day,y=ub_pwlrdx_3] {\data};
        \end{axis}
    \end{tikzpicture}
    \begin{tikzpicture}[baseline=(current bounding box.north)]
        \begin{axis}[enlargelimits=false,height=4.3cm,xlabel=Price Profile (Day),ymajorticks=false,
            xticklabels={1,2,3,4,5},xtick={1,2,3,4,5},width=0.37\linewidth,xtick pos=left,
            ytick pos=left,xmin=1,xmax=5,ymin=0.0,ymax=2.25]
            \node[anchor=north east] at (rel axis cs:0.99,0.98) {$K = 24$};
            \pgfplotstableread[col sep = comma]{MATM-24_Steps-UB.csv}{\data};
            \addplot[index of colormap=0 of Dark2-8,mark=*,only marks,mark size=3pt,opacity=0.5] table [x=day,y=ub_lrdx_1] {\data};
            \addplot[index of colormap=1 of Dark2-8,mark=*,only marks,mark size=3pt,opacity=0.5] table [x=day,y=ub_lrdx_2] {\data};
            \addplot[index of colormap=2 of Dark2-8,mark=*,only marks,mark size=3pt,opacity=0.5] table [x=day,y=ub_lrdx_3] {\data};
            \addplot[index of colormap=0 of Dark2-8,mark=diamond*,only marks,mark size=3pt,opacity=0.5] table [x=day,y=ub_pwlrdx_1] {\data};
            \addplot[index of colormap=1 of Dark2-8,mark=diamond*,only marks,mark size=3pt,opacity=0.5] table [x=day,y=ub_pwlrdx_2] {\data};
            \addplot[index of colormap=2 of Dark2-8,mark=diamond*,only marks,mark size=3pt,opacity=0.5] table [x=day,y=ub_pwlrdx_3] {\data};
        \end{axis}
    \end{tikzpicture}
    \begin{tikzpicture}[baseline=(current bounding box.north)]
        \begin{axis}[enlargelimits=false,height=4.3cm,ymajorticks=false,
            xlabel=Price Profile (Day),xticklabels={1,2,3,4,5},xtick={1,2,3,4,5},
            width=0.37\linewidth,xtick pos=left,ytick pos=left,xmin=1,xmax=5,ymin=0.0,ymax=2.25]
            \node[anchor=north east] at (rel axis cs:0.99,0.98) {$K = 48$};
            \pgfplotstableread[col sep = comma]{MATM-48_Steps-UB.csv}{\data};
            \addplot[index of colormap=0 of Dark2-8,mark=*,only marks,mark size=3pt,opacity=0.5] table [x=day,y=ub_lrdx_1] {\data};
            \addplot[index of colormap=1 of Dark2-8,mark=*,only marks,mark size=3pt,opacity=0.5] table [x=day,y=ub_lrdx_2] {\data};
            \addplot[index of colormap=2 of Dark2-8,mark=*,only marks,mark size=3pt,opacity=0.5] table [x=day,y=ub_lrdx_3] {\data};
            \addplot[index of colormap=0 of Dark2-8,mark=diamond*,only marks,mark size=3pt,opacity=0.5] table [x=day,y=ub_pwlrdx_1] {\data};
            \addplot[index of colormap=1 of Dark2-8,mark=diamond*,only marks,mark size=3pt,opacity=0.5] table [x=day,y=ub_pwlrdx_2] {\data};
            \addplot[index of colormap=2 of Dark2-8,mark=diamond*,only marks,mark size=3pt,opacity=0.5] table [x=day,y=ub_pwlrdx_3] {\data};
        \end{axis}
    \end{tikzpicture}
    \caption{Comparison of primal bound quality based on relaxation model and resolution, $\xi$, for \texttt{AT(M)} instances.}
    \label{figure:primal-bound-atm}%
\end{figure}

Figures \ref{figure:primal-bound-atm} and \ref{figure:primal-bound-poormond} compare the effectiveness of the six relaxations in producing \eqref{equation:mincp}-feasible primal bounds for \texttt{AT(M)} and \texttt{Poormond} instances.
A figure comparing \texttt{Simple FSD} solutions is not included due to the similar convergence properties observed across all relaxations.
First, inspecting \texttt{AT(M)} solutions for $K = 12$, nearly all methods produce the same upper bounds over all instances, although \eqref{equation:milp-pw} relaxations with $\xi = 1$ meter do not find solutions for any instances.
For $K = 24$, \eqref{equation:milp-oa}-based relaxations are more capable of producing good solutions when compared with \eqref{equation:milp-pw}.
For $K = 48$, \eqref{equation:milp-pw} relaxations appear to be too challenging for discovering any feasible solution, whereas \eqref{equation:milp-oa}-based relaxations consistently find solutions within $2\%$ of one another.

\begin{figure}[t]
    \centering
    \begin{tikzpicture}[baseline=(current bounding box.north)]
        \begin{axis}[ylabel={$\%$ from Best},enlargelimits=false,height=4.3cm,
            legend style={at={(0.05, 0.94)},anchor=north west,font=\footnotesize},
            xlabel=Price Profile (Day),xticklabels={1,2,3,4,5},xtick={1,2,3,4,5},
            width=0.37\linewidth,xtick pos=left,ytick pos=left,xmin=1,xmax=5,
            ymin=0.0,ymax=5.0,legend columns=1,legend entries={$\xi = 1$ m, $\xi = 5$ m,
            $\xi = 25$ m,\eqref{equation:milp-oa},\eqref{equation:milp-pw}},legend cell align={left}]
            \node[anchor=north east] at (rel axis cs:0.99,0.98) {$K = 12$};
            \addlegendimage{index of colormap=0 of Dark2-8,mark=square*,only marks,opacity=0.5};
            \addlegendimage{index of colormap=1 of Dark2-8,mark=square*,only marks,opacity=0.5};
            \addlegendimage{index of colormap=2 of Dark2-8,mark=square*,only marks,opacity=0.5};
            \addlegendimage{black,mark=*,only marks};
            \addlegendimage{black,mark=diamond*,only marks};
            \pgfplotstableread[col sep = comma]{Poormond-12_Steps-UB.csv}{\data};
            \addplot[index of colormap=0 of Dark2-8,mark=*,only marks,mark size=3pt,opacity=0.5] table [x=day,y=ub_lrdx_1] {\data};
            \addplot[index of colormap=1 of Dark2-8,mark=*,only marks,mark size=3pt,opacity=0.5] table [x=day,y=ub_lrdx_2] {\data};
            \addplot[index of colormap=2 of Dark2-8,mark=*,only marks,mark size=3pt,opacity=0.5] table [x=day,y=ub_lrdx_3] {\data};
            \addplot[index of colormap=0 of Dark2-8,mark=diamond*,only marks,mark size=3pt,opacity=0.5] table [x=day,y=ub_pwlrdx_1] {\data};
            \addplot[index of colormap=1 of Dark2-8,mark=diamond*,only marks,mark size=3pt,opacity=0.5] table [x=day,y=ub_pwlrdx_2] {\data};
            \addplot[index of colormap=2 of Dark2-8,mark=diamond*,only marks,mark size=3pt,opacity=0.5] table [x=day,y=ub_pwlrdx_3] {\data};
        \end{axis}
    \end{tikzpicture}
    \begin{tikzpicture}[baseline=(current bounding box.north)]
        \begin{axis}[enlargelimits=false,height=4.3cm,xlabel=Price Profile (Day),ymajorticks=false,
            xticklabels={1,2,3,4,5},xtick={1,2,3,4,5},width=0.37\linewidth,xtick pos=left,
            ytick pos=left,xmin=1,xmax=5,ymin=0.0,ymax=5.0]
            \node[anchor=north east] at (rel axis cs:0.99,0.98) {$K = 24$};
            \pgfplotstableread[col sep = comma]{Poormond-24_Steps-UB.csv}{\data};
            \addplot[index of colormap=0 of Dark2-8,mark=*,only marks,mark size=3pt,opacity=0.5] table [x=day,y=ub_lrdx_1] {\data};
            \addplot[index of colormap=1 of Dark2-8,mark=*,only marks,mark size=3pt,opacity=0.5] table [x=day,y=ub_lrdx_2] {\data};
            \addplot[index of colormap=2 of Dark2-8,mark=*,only marks,mark size=3pt,opacity=0.5] table [x=day,y=ub_lrdx_3] {\data};
            \addplot[index of colormap=0 of Dark2-8,mark=diamond*,only marks,mark size=3pt,opacity=0.5] table [x=day,y=ub_pwlrdx_1] {\data};
            \addplot[index of colormap=1 of Dark2-8,mark=diamond*,only marks,mark size=3pt,opacity=0.5] table [x=day,y=ub_pwlrdx_2] {\data};
            \addplot[index of colormap=2 of Dark2-8,mark=diamond*,only marks,mark size=3pt,opacity=0.5] table [x=day,y=ub_pwlrdx_3] {\data};
        \end{axis}
    \end{tikzpicture}
    \begin{tikzpicture}[baseline=(current bounding box.north)]
        \begin{axis}[enlargelimits=false,height=4.3cm,ymajorticks=false,
            xlabel=Price Profile (Day),xticklabels={1,2,3,4,5},xtick={1,2,3,4,5},
            width=0.37\linewidth,xtick pos=left,ytick pos=left,xmin=1,xmax=5,ymin=0.0,ymax=5.0]
            \node[anchor=north east] at (rel axis cs:0.99,0.98) {$K = 48$};
            \pgfplotstableread[col sep = comma]{Poormond-48_Steps-UB.csv}{\data};
            \addplot[index of colormap=0 of Dark2-8,mark=*,only marks,mark size=3pt,opacity=0.5] table [x=day,y=ub_lrdx_1] {\data};
            \addplot[index of colormap=1 of Dark2-8,mark=*,only marks,mark size=3pt,opacity=0.5] table [x=day,y=ub_lrdx_2] {\data};
            \addplot[index of colormap=2 of Dark2-8,mark=*,only marks,mark size=3pt,opacity=0.5] table [x=day,y=ub_lrdx_3] {\data};
            \addplot[index of colormap=0 of Dark2-8,mark=diamond*,only marks,mark size=3pt,opacity=0.5] table [x=day,y=ub_pwlrdx_1] {\data};
            \addplot[index of colormap=1 of Dark2-8,mark=diamond*,only marks,mark size=3pt,opacity=0.5] table [x=day,y=ub_pwlrdx_2] {\data};
            \addplot[index of colormap=2 of Dark2-8,mark=diamond*,only marks,mark size=3pt,opacity=0.5] table [x=day,y=ub_pwlrdx_3] {\data};
        \end{axis}
    \end{tikzpicture}
    \caption{Comparison of primal bound quality based on relaxation model and resolution, $\xi$, for \texttt{Poormond} instances.}
    \label{figure:primal-bound-poormond}%
\end{figure}

In Figure \ref{figure:primal-bound-poormond}, for $K = 12$, all relaxations find primal solutions for all instances, although tighter relaxations sometimes produce better solutions.
This is apparent in price profile $1$, where tight \eqref{equation:milp-pw} relaxations outperform \eqref{equation:milp-oa} relaxations.
For $K = 24$, all relaxations produce good primal solutions within the alotted time limit.
For $K = 48$, tight \eqref{equation:milp-pw} relaxations are generally more capable of finding better upper bounds, although \eqref{equation:milp-oa}-based relaxations also consistently find good solutions.
Figures \ref{figure:primal-bound-atm} and \ref{figure:primal-bound-poormond} suggest there is an important tradeoff between the complexity of a MILP model and the ability to explore integer solutions.
In general, \eqref{equation:milp-oa} relaxations with $\xi \in \{1, 5\}$ and continuously relaxed direction variables produce good feasible solutions for all instances.

\begin{table}[t]
    \small
    \centering
    \begin{center}
    \begin{tabular}{cc|cccc|cccc|cccc|}
        \cline{3-14} & & \multicolumn{4}{c|}{$K = 12$} &
                         \multicolumn{4}{c|}{$K = 24$} &
                         \multicolumn{4}{c|}{$K = 48$} \\
        \cline{1-14}
        & Day & $\textrm{UB}$ & $\textrm{LB}$ & Gap & Time & $\textrm{UB}$ & $\textrm{LB}$ & Gap & Time & $\textrm{UB}$ & $\textrm{LB}$ & Gap & Time \\
        \cline{1-14}
        \parbox[t]{2mm}{\multirow{5}{*}{\rotatebox[origin=c]{90}{\texttt{Simple FSD}}}}
        & \multicolumn{1}{|c|}{1} & - & - & - & - & 155.6 & 155.6 & 0.0\% & 0.5 s & 152.9 & 152.8 & 0.1\% & 5.5 s \\
        & \multicolumn{1}{|c|}{2} & - & - & - & - & 159.1 & 159.0 & 0.0\% & 0.4 s & 157.0 & 156.9 & 0.0\% & 5.8 s \\
        & \multicolumn{1}{|c|}{3} & - & - & - & - & 172.7 & 172.6 & 0.0\% & 0.5 s & 170.6 & 170.4 & 0.1\% & 6.5 s \\
        & \multicolumn{1}{|c|}{4} & - & - & - & - & 182.2 & 182.1 & 0.0\% & 0.4 s & 178.0 & 177.8 & 0.1\% & 2.2 s \\
        & \multicolumn{1}{|c|}{5} & - & - & - & - & 149.1 & 149.0 & 0.1\% & 0.5 s & 147.3 & 147.1 & 0.1\% & 2.5 s \\
        \cline{1-14}
        \parbox[t]{2mm}{\multirow{5}{*}{\rotatebox[origin=c]{90}{\texttt{AT(M)}}}}
        & \multicolumn{1}{|c|}{1} & 816.9 & 779.6 & 4.6\% & 89.7 s & 742.4 & 729.9 & 1.7\% & 2177.2 s & 738.1 & 719.3 & 2.5\% & Lim. \\
        & \multicolumn{1}{|c|}{2} & 821.4 & 782.5 & 4.7\% & 83.0 s & 737.1 & 721.9 & 2.1\% & 1874.3 s & 734.5 & 707.8 & 3.6\% & Lim. \\
        & \multicolumn{1}{|c|}{3} & 853.2 & 813.3 & 4.7\% & 75.1 s & 769.8 & 753.8 & 2.1\% & 1376.3 s & 764.1 & 738.4 & 3.4\% & Lim. \\
        & \multicolumn{1}{|c|}{4} & 947.6 & 902.3 & 4.8\% & 85.3 s & 829.0 & 814.1 & 1.8\% & 2011.8 s & 821.9 & 795.5 & 3.2\% & Lim. \\
        & \multicolumn{1}{|c|}{5} & 988.0 & 938.2 & 5.0\% & 69.0 s & 694.2 & 681.3 & 1.9\% & 1969.2 s & 692.8 & 651.9 & 5.9\% & Lim. \\
        \cline{1-14}
        \parbox[t]{2mm}{\multirow{5}{*}{\rotatebox[origin=c]{90}{\texttt{Poormond}}}}
        & \multicolumn{1}{|c|}{1} & 119.1 & 115.3 & 3.2\% & Lim. & 114.5 & 110.1 & 3.8\% & Lim. & 112.4 & 109.0 & 3.0\% & Lim. \\
        & \multicolumn{1}{|c|}{2} & 120.5 & 118.5 & 1.7\% & Lim. & 117.1 & 112.4 & 4.1\% & Lim. & 116.5 & 111.2 & 4.6\% & Lim. \\
        & \multicolumn{1}{|c|}{3} & 133.6 & 130.8 & 2.1\% & Lim. & 127.8 & 124.3 & 2.8\% & Lim. & 127.1 & 122.6 & 3.5\% & Lim. \\
        & \multicolumn{1}{|c|}{4} & 147.0 & 143.2 & 2.6\% & Lim. & 140.6 & 137.8 & 2.0\% & Lim. & 140.7 & 135.9 & 3.4\% & Lim. \\
        & \multicolumn{1}{|c|}{5} & 131.4 & 126.4 & 3.9\% & Lim. & 100.2 & 98.2 & 2.0\% & Lim. & 102.8 & 92.5 & 10.0\% & Lim. \\
        \cline{1-14}
    \end{tabular}
    \end{center}
    \caption{Primal and dual solution statistics after one hour when using \eqref{equation:milp-oa} with $\xi = 1$ meter.}
    \label{table:formulation_upper_bound_improvements}
\end{table}

In Table \ref{table:formulation_upper_bound_improvements}, solution statistics are shown for all instances solved using \eqref{equation:milp-oa} with $\xi = 1$ meter and continuously relaxed direction variables.
First, for \texttt{Simple FSD}, all $K = 12$ instances are proven infeasible by the relaxation, as indicated by the dashed entries.
For the remaining \texttt{Simple FSD} instances (i.e., $K \in \{24, 48\}$), each instance is solved within seven seconds.
We emphasize that, here, each upper bound (``UB'') is a \emph{true upper bound} of \eqref{equation:mincp}, and each lower bound (``LB'') is a \emph{lower bound} of \eqref{equation:milp-oa}, and thus also \eqref{equation:mincp}.
However, because termination depends only upon convergence of \eqref{equation:milp-oa}, optimal lower bounds do not necessarily coincide with \eqref{equation:mincp}-feasible upper bounds.

For \texttt{AT(M)}, all $K = 12$ instances converge within ninety seconds, and all $K = 24$ instances converge within the one hour time limit.
For $K = 48$, none of the relaxations converge, although optimality gaps reported at termination are generally small and less than $6\%$.
For \texttt{Poormond}, even for the smallest instances ($K = 12$), none of the relaxations converge, although reported optimality gaps at termination are typically a few percent or less.
The exception is the ``Day 5'' instance for $K = 48$, which reports a $10\%$ gap.
It is nonetheless remarkable that (i) \emph{\eqref{equation:mincp}-feasible solutions are discovered for all instances within the one hour time limit without a specialized heuristic}, and (ii) \emph{reported optimality gaps are generally small, even given the small time limit and relative looseness of the relaxation}.

To further evaluate our modeling enhancements, we attempt to draw some comparisons with results presented by \citet{bonvin2021pump}.
We first remark that \citet{bonvin2021pump} generally employ an upper-bounding heuristic that allows the discreteness of control decisions to be violated.
This makes a direct comparison challenging.
However, they also provide one tabulation of results for \texttt{Poormond} instances where discreteness \emph{is} required.
Here, for an equal time limit, they find feasible solutions for only nine of fifteen instances, whereas we discover feasible solutions for all instances.
For $K = 12$, they report an average optimality gap of $7.5\%$ across instances where upper bounds are discovered, whereas we report an average gap of $2.7\%$.
For $K = 24$, their average gap is $9.0\%$, whereas ours is $2.9\%$.
Finally, for $K = 48$, their average is $8.5\%$, whereas ours is $4.9\%$.
These comparisons imply the techniques described in our paper may provide tighter relaxations that (i) are more capable of producing feasible solutions and (ii) substantially improve OWF dual bounds.

\subsection{Relative Effects of Formulation Improvements}
\label{subsection:relative-effects}
Sections \ref{subsection:effects_of_obbt} through \ref{subsection:effects_of_valid_inequalities} evaluated the effects of bound tightening and valid inequalities on the OWF's root relaxation, while Section \ref{subsection:effects_of_polyhedral_parameterizations} evaluated the ability of various polyhedral relaxations to produce primal solutions that are feasible to the original OWF.
These sections do not, however, evaluate the effects of bound tightening and valid inequalities on the overall convergence of the LP/NLP-BB algorithm used in Section \ref{subsection:effects_of_polyhedral_parameterizations}.
In this section, we evaluate the effects of each model strengthening technique sequentially on the five \texttt{Poormond}, $K = 48$ instances, which are the most challenging instances of this study.

\begin{table}[t]
    \footnotesize
    \centering
    \begin{center}
    \begin{tabular}{c|cc|cc|cc|cc|cc|cc|cc|cc|}
        \cline{2-17} & \multicolumn{2}{c|}{Naive} &
                       \multicolumn{2}{c|}{BT-SR} & \multicolumn{2}{c|}{BT-SS} &
                       \multicolumn{2}{c|}{BT-SQ} & \multicolumn{2}{c|}{BT-TR} &
                       \multicolumn{2}{c|}{BT-TS} & \multicolumn{2}{c|}{Dual Cuts} & \multicolumn{2}{c|}{OBCG} \\
        \cline{1-17}
        \multicolumn{1}{|c|}{Day} & LB & UB & LB & UB & LB & UB & LB & UB & LB & UB & LB & UB & LB & UB & LB & UB \\
        \cline{1-17}
        \multicolumn{1}{|c|}{1} & $103.4$ & $-$ & $103.7$ & $-$ & $103.9$ & $-$ & $103.6$ & $116.0$ & $103.4$ & $113.9$ & $103.7$ & $116.1$ & $105.0$ & $-$ & $109.0$ & $112.4$ \\
        \multicolumn{1}{|c|}{2} & $105.5$ & $-$ & $105.9$ & $-$ & $105.6$ & $-$ & $105.8$ & $117.0$ & $105.8$ & $116.4$ & $105.7$ & $-$ & $107.4$ & $-$ & $111.2$ & $116.5$ \\
        \multicolumn{1}{|c|}{3} & $116.6$ & $-$ & $117.1$ & $-$ & $117.0$ & $-$ & $117.5$ & $129.3$ & $117.9$ & $126.8$ & $117.9$ & $126.5$ & $119.0$ & $-$ & $122.6$ & $127.1$ \\
        \multicolumn{1}{|c|}{4} & $128.6$ & $-$ & $129.5$ & $-$ & $129.5$ & $-$ & $129.7$ & $143.5$ & $129.6$ & $142.9$ & $129.0$ & $146.5$ & $131.2$ & $141.9$ & $135.9$ & $140.7$ \\
        \multicolumn{1}{|c|}{5} & $89.1$ & $-$ & $89.3$ & $-$ & $89.2$ & $-$ & $89.3$ & $103.4$ & $89.7$ & $104.4$ & $89.7$ & $102.7$ & $89.6$ & $107.5$ & $92.5$ & $102.8$ \\
        \cline{1-17}
    \end{tabular}
    \end{center}
    \caption{Evaluation of \texttt{Poormond} ($K = 48$) primal and dual bounds for \eqref{equation:milp-oa} with $\xi = 1$ meter.}
    \label{table:bound-improvements-poormond}
\end{table}

Table \ref{table:bound-improvements-poormond} summarizes the lower and upper bounds achieved using each model strengthening technique with the \eqref{equation:milp-oa} relaxation studied in Section \ref{subsection:effects_of_polyhedral_parameterizations}, where discrete direction variables are continuously relaxed and $\xi = 1$ meter.
Each header column (e.g., BT-SR, BT-SS) assumes use of the bound tightening or cutting plane technique that precedes it (e.g., ``Naive'' and BT-SR, respectively).
It is first apparent that the effects of bound tightening are often marginal until BT-SQ is used, which enables the discovery of many feasible primal solutions.
The effects of BT-TR and BT-TS on the lower bound appear small and mostly inconsistent.
Duality-based cuts generally provide lower bound improvements around $1\%$.
Finally, application of OBCG cuts provides relatively large lower bound improvements, i.e., $3\%$ to $4\%$, and these cuts also enable the discovery of consistently good feasible solutions.

\begin{figure}[t]
    \centering
    \begin{tikzpicture}[baseline=(current bounding box.north)]
        \begin{semilogxaxis}[ylabel={Dual (lower) bound},enlargelimits=false,height=6.0cm,
            legend style={at={(0.03, 0.82)},anchor=north west,font=\tiny},
            xlabel=Wall clock time (s),width=0.49\linewidth,xtick pos=left,ytick pos=left,
            legend columns=3,legend cell align={left}]
            \pgfplotstableread[col sep = comma]{Poormond-48_Steps-Day_24-UB-NO-BT.csv}{\data};
            \addplot[thick,black,dashed] table [x=time_elapsed,y=lower_bound]{\data};
            \addlegendentry{Naive};

            \pgfplotstableread[col sep = comma]{Poormond-48_Steps-Day_24-UB-BT-SS-R.csv}{\data};
            \addplot[thick,index of colormap=0 of Dark2-8] table [x=time_elapsed,y=lower_bound]{\data};
            \addlegendentry{BT-SR};

            \pgfplotstableread[col sep = comma]{Poormond-48_Steps-Day_24-UB-BT-SS.csv}{\data};
            \addplot[thick,index of colormap=1 of Dark2-8] table [x=time_elapsed,y=lower_bound]{\data};
            \addlegendentry{BT-SS};

            \pgfplotstableread[col sep = comma]{Poormond-48_Steps-Day_24-UB-BT-SQ.csv}{\data};
            \addplot[thick,index of colormap=2 of Dark2-8] table [x=time_elapsed,y=lower_bound]{\data};
            \addlegendentry{BT-SQ};

            \pgfplotstableread[col sep = comma]{Poormond-48_Steps-Day_24-UB-BT-TS-R.csv}{\data};
            \addplot[thick,index of colormap=3 of Dark2-8] table [x=time_elapsed,y=lower_bound]{\data};
            \addlegendentry{BT-TR};

            \pgfplotstableread[col sep = comma]{Poormond-48_Steps-Day_24-UB-BT-TS.csv}{\data};
            \addplot[thick,index of colormap=4 of Dark2-8] table [x=time_elapsed,y=lower_bound]{\data};
            \addlegendentry{BT-TS};

            \pgfplotstableread[col sep = comma]{Poormond-48_Steps-Day_24-UB-DUAL-CUTS.csv}{\data};
            \addplot[thick,index of colormap=5 of Dark2-8] table [x=time_elapsed,y=lower_bound]{\data};
            \addlegendentry{Dual Cuts};

            \pgfplotstableread[col sep = comma]{Poormond-48_Steps-Day_24-UB-ALL-CUTS.csv}{\data};
            \addplot[thick,index of colormap=6 of Dark2-8] table [x=time_elapsed,y=lower_bound]{\data};
            \addlegendentry{OBCG};
        \end{semilogxaxis}
    \end{tikzpicture}
    \begin{tikzpicture}[baseline=(current bounding box.north)]
        \begin{axis}[enlargelimits=false,height=6.0cm,
            legend style={at={(0.05, 0.94)},ylabel=Primal (upper) bound,anchor=north west,font=\footnotesize},
            xlabel=Wall clock time (s),width=0.49\linewidth,xtick pos=left,ytick pos=left,
            legend columns=1,legend cell align={left},xmin=1.0,ymin=139]
            \pgfplotstableread[col sep = comma]{Poormond-48_Steps-Day_24-UB-NO-BT.csv}{\data};
            \addplot[thick,black,dashed] table [x=time_elapsed,y=upper_bound]{\data};

            \pgfplotstableread[col sep = comma]{Poormond-48_Steps-Day_24-UB-BT-SS-R.csv}{\data};
            \addplot[thick,index of colormap=0 of Dark2-8] table [x=time_elapsed,y=upper_bound]{\data};

            \pgfplotstableread[col sep = comma]{Poormond-48_Steps-Day_24-UB-BT-SS.csv}{\data};
            \addplot[thick,index of colormap=1 of Dark2-8] table [x=time_elapsed,y=upper_bound]{\data};

            \pgfplotstableread[col sep = comma]{Poormond-48_Steps-Day_24-UB-BT-SQ.csv}{\data};
            \addplot[thick,index of colormap=2 of Dark2-8] table [x=time_elapsed,y=upper_bound]{\data};

            \pgfplotstableread[col sep = comma]{Poormond-48_Steps-Day_24-UB-BT-TS-R.csv}{\data};
            \addplot[thick,index of colormap=3 of Dark2-8] table [x=time_elapsed,y=upper_bound]{\data};

            \pgfplotstableread[col sep = comma]{Poormond-48_Steps-Day_24-UB-BT-TS.csv}{\data};
            \addplot[thick,index of colormap=4 of Dark2-8] table [x=time_elapsed,y=upper_bound]{\data};

            \pgfplotstableread[col sep = comma]{Poormond-48_Steps-Day_24-UB-DUAL-CUTS.csv}{\data};
            \addplot[thick,index of colormap=5 of Dark2-8] table [x=time_elapsed,y=upper_bound]{\data};

            \pgfplotstableread[col sep = comma]{Poormond-48_Steps-Day_24-UB-ALL-CUTS.csv}{\data};
            \addplot[thick,index of colormap=6 of Dark2-8] table [x=time_elapsed,y=upper_bound]{\data};
        \end{axis}
    \end{tikzpicture}
    \caption{Dual and primal bound convergence profiles, respectively, for the Poormond, $K = 48$, ``Day 4'' instance.}
    \label{figure:convergence-poormond}%
\end{figure}
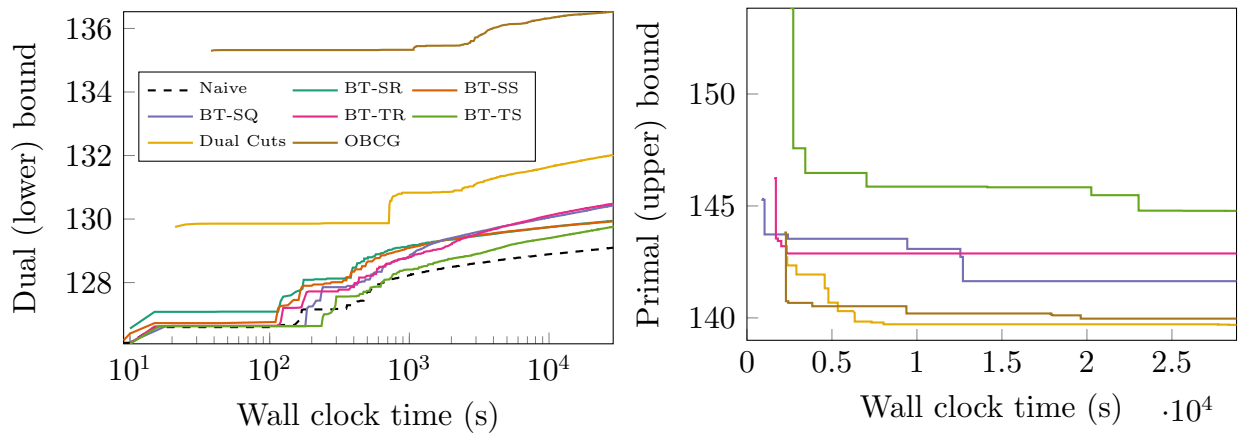

Figure \ref{figure:convergence-poormond} illustrates dual and primal bound convergence, respectively, when applying each method to the ``Day 4'' instance of Table \ref{table:formulation_upper_bound_improvements}, using a wall clock time limit of eight hours (instead of one hour).
First, over time, dual bounds tend to become larger as more strengthening techniques are applied.
More intensive OBBT methods (e.g., BT-SQ and BT-TR) eventually show moderate improvements over other methods (i.e., BT-SR and BT-SS).
Further, the dual bound improvements are the most dramatic after applying duality- and OBCG-based cuts.
For primal bounds, production of feasible solutions begins after BT-SQ is used, although the quality of incumbent solutions that are discovered tends to vary.
Nonetheless, it is interesting to observe that as more strengthening techniques are applied (e.g., the dual cuts and OBCG lines), the quality of primal solutions is improved.

It is challenging to draw conclusions about the value of each technique in isolation, as tighter variable bounds can entail stronger cuts, and the combination of duality-based and OBCG cuts may provide more value than if they were applied independently.
Furthermore, improvements discovered by \textsc{Gurobi} during the search (e.g., via presolve, heuristics, and cutting planes) are difficult to anticipate a priori.
Nonetheless, using Table \ref{table:bound-improvements-poormond} and Figure \ref{figure:convergence-poormond} as examples, we summarize that (i) each technique helps improve the lower bound, with OBCG providing the largest lower bound improvements, (ii) application of these model strengthening techniques leads to the discovery of higher-quality primal solutions, and (iii) it is the combination of all techniques that contributes to the results presented here.

\section{Conclusion}
\label{section:conclusion}
This paper explored relaxation-based solution techniques for the OWF problem, which minimizes the cost of pump energy consumption over a fixed time horizon.
Unlike other studies, this paper focused primarily on formulation, preprocessing, and algorithmic techniques that \emph{strengthen relaxations} of the OWF, as slow dual bound improvement has been identified as a weakness of recent optimization methods (e.g., \citealp{vieira2020optimizing,bonvin2021pump}).
First, the OWF was formulated as a MINLP.
To better leverage modern MILP solvers, two polyhedral relaxation-based formulations were then developed, each of whose accuracy and tractability are controlled by the partitioning of flow values used to model nonlinearities.
These contributions strengthen formulations used in the previous literature.

Expanding upon the previous literature, novel duality-based convex valid inequalities, similar to those derived for WDN design by \cite{tasseff2020exact}, were developed for use in the OWF problem.
Additionally, a number of preprocessing techniques were developed to increase the tightness of relaxed OWF problem specifications.
First, multiple novel OBBT procedures were introduced, and their efficacies and computational tradeoffs were evaluated.
Then, novel OBCG methods were developed that aim at further strengthening OWF relaxations.
Empirical results indicate that our model strengthening techniques can substantially improve the tightness of an OWF relaxation, especially for more complex instances involving many unique controllable components (e.g., pumps) and storage tanks.

Finally, a variety of polyhedral relaxations were evaluated to explore the computational tradeoffs between OWF model accuracy and the ability to efficienty discover feasible primal solutions.
These results indicate that (i) the best relaxation parameterization often depends on instance complexity, network topology, and the tightness of variable bounds; (ii) highly accurate models are sometimes prohibitively expensive; and (iii) exploring integer solutions (e.g., by using a less accurate relaxation) is often crucial for efficiently discovering feasible primal solutions.
Comparing our results with the prior state of the art \citep{bonvin2021pump} suggests that the model strengthening techniques we propose, as well as a careful polyhedral partitioning of nonlinearities, provide relaxations that are capable of finding good primal solutions for complex OWFs while also significantly reducing optimality gaps.

We suggest that future work should focus on the development of heuristics that generate solutions to \eqref{equation:mincp}.
The availability of even one primal solution could further constrain OBBT and OBCG subproblems, which can lead to substantial improvements in variable bounds and strengthening of cuts.
As mentioned in Section \ref{section:literature_review}, however, many existing techniques rely on metaheuristics, which allow infeasibilities, or mathematical programming methods, which rely on specific assumptions concerning WDN topology.
One promising general mathematical programming heuristic has recently been proposed by \cite{bonvin2019extended}, in which low-cost feasible solutions are obtained by solving a relaxed problem with preprocessed exact surrogate data that assume particular tank levels.
Nonetheless, success of this approach appears highly dependent on the selection of fixed tank levels and variations in OWF demand and price profiles.
We strongly suggest that new heuristics should ensure that (i) solutions are feasible to \eqref{equation:mincp} and (ii) they are generally applicable to a variety of WDN topologies and OWF parameterizations.

\ACKNOWLEDGMENT{%
This work is supported by the U.S. Department of Energy's (DOE's) Advanced Grid Modeling (AGM) project \emph{Coordinated Planning and Operation of Water and Power Infrastructures for Increased Resilience and Reliability}.
Incorporation of the \textsc{PolyhedralRelaxations} \textsc{Julia} package was supported by Los Alamos National Laboratory's Directed Research and Development program under the project \emph{Fast, Linear Programming-Based Algorithms with Solution Quality Guarantees for Nonlinear Optimal Control Problems} (20220006ER).
All work at Los Alamos National Laboratory was conducted under the auspices of the National Nuclear Security Administration of the U.S. DOE under Contract No. 89233218CNA000001. This work was authored in part by the National Renewable Energy Laboratory, operated by Alliance for Sustainable Energy, LLC, for the U.S. DOE under Contract No. DE-AC36-08GO28308. 
We gratefully acknowledge Marina A. Epelman at the University of Michigan for her thoughtful suggestions regarding an early draft of this paper.
We also thank our AGM program manager, Alireza Ghassemian, for his support of this research.

The views expressed in the article do not necessarily represent the views of the DOE or the U.S. Government.
The U.S. Government retains and the publisher, by accepting the article for publication, acknowledges that the U.S. Government retains a nonexclusive, paid-up, irrevocable, worldwide license to publish or reproduce the published form of this work, or allow others to do so, for U.S. Government purposes. %
}

\bibliographystyle{informs2014}
\bibliography{bibliography.bib}

\clearpage

\begin{APPENDICES}
\section{Direction-based Valid Inequalities}
Here, we detail the direction-based valid inequalities briefly mentioned near the end of Section 4.1.
These inequalities are used to strengthen direction-based OWF formulations.
We note that they are similar to the inequalities used by \cite{borraz2016convex} for natural gas network modeling.
The first constraints are
\begin{equation}
    \sum_{\mathclap{(i, j) \in \delta_{i}^{+}}} y_{ij}^{k} \geq 1, \, \forall i \in \mathcal{R} \cup \{i^{\prime} \in \mathcal{D} : \overline{q}_{i^{\prime}}^{k} > 0\}, \, \forall k \in \mathcal{K} \label{eqn:source-flow},
\end{equation}
which model that at least one link must send water \emph{away} from each source in the WDN.
Similarly, the following valid inequalities are added for every node in the WDN that \emph{withdraws} water from the network:
\begin{equation}
    \sum_{\mathclap{(j, i) \in \delta_{i}^{-}}} y_{ji}^{k} + \sum_{\mathclap{(i, j) \in \delta_{i}^{+}}} (1 - y_{ij}^{k}) \geq 1, \, \forall i \in \mathcal{D}, \, \forall k \in \mathcal{K} : \overline{q}_{i}^{k} < 0 \label{eqn:demand-flow}.
\end{equation}
These constraints model that at least one link must provide water \emph{to} each consumption node.
Finally,
\begin{subequations}\begin{align}
    \sum_{\mathclap{(j, i) \in \delta_{i}^{-}}} y_{ji}^{k} - \sum_{\mathclap{(i, j) \in \delta_{i}^{+}}} y_{ij}^{k} &= 0, \, \forall i \in \mathcal{D}, \, \forall k \in \mathcal{K} : \left(\overline{q}_{i}^{k} = 0\right) \land \left(\deg^{+}_{i} = \deg^{-}_{i} = 1\right) \label{eqn:deg-2-flow-1} \\
    \sum_{\mathclap{(j, i) \in \delta_{i}^{-}}} y_{ji}^{k} + \sum_{\mathclap{(i, j) \in \delta_{i}^{+}}} y_{ij}^{k} &= 1, \, \forall i \in \mathcal{D}, \, \forall k \in \mathcal{K} : \left(\overline{q}_{i}^{k} = 0\right) \land \left(\deg^{\pm}_{i} = 2\right) \land \left(\deg^{\mp}_{i} = 0\right) \label{eqn:deg-2-flow-2}
\end{align}\label{eqn:deg-2-flow}\end{subequations}
model directionality at zero demand junction nodes with a degree of exactly two.
The implication is that, for a pass-through junction, the direction of incoming flow must be equal to the direction of outgoing flow.

\section{The Content Model}
Section \ref{section:feasibility-step} states a set of novel valid inequalities derived using an approach that reformulates nonlinear physical constraints.
In this and subsequent appendices, we provide details regarding our derivation for the interested reader.
These details begin with a discussion of the known ``Content Model'' for WDN \emph{analysis}.

In the context of network analysis (e.g., the simulation of a network with fixed control decisions), a relaxed version of \eqref{equation:mincp} can be solved via an explicit sequential method while neglecting flow and head bounds.
At the core of such methods is a technique that solves a convex programming problem and provides a solution satisfying Constraints \eqref{equation:mincp-pipe-head-loss}, \eqref{equation:mincp-pump-head-gain}, \eqref{equation:mincp-pump-head}, and \eqref{equation:flow-conservation} (i.e., nonlinear constraints) for one $k \in \mathcal{K}$.
For pipe-only networks, these convex ``Content'' and ``Co-Content'' models originated with \cite{cherry1951cxvii} and \cite{collins1978solving} and were recently employed by \cite{tasseff2020exact} for optimal gravity-fed pipe network design.

Here, we explicitly extend the Content Model \citep{collins1978solving} to include pumps and valves as follows:
\begin{equation}
\tag{P$^{k}$}
\begin{aligned}
    & \text{minimize}
    & &
    \begin{gathered}
        \sum_{(i, j) \in \mathcal{A}} \frac{L_{ij} r_{ij}}{1 + \alpha} \left[(q_{ij}^{k +})^{1 + \alpha} + (q_{ij}^{k -})^{1 + \alpha}\right]
        - \sum_{(i, j) \in \mathcal{P}^{k}_{1}} \left[a_{ij} q_{ij}^{k +} + \frac{b_{ij} (q_{ij}^{k +})^{c_{ij} + 1}}{c_{ij} + 1} \right] \\
        - \sum_{i \in \mathcal{R} \cup \mathcal{T}} h_{i}^{k} \left[\sum_{(i, j) \in \delta^{k +}_{i}} (q_{ij}^{k +} - q_{ij}^{k -}) - \sum_{(j, i) \in \delta^{k -}_{i}} (q_{ji}^{k +} - q_{ji}^{k -})\right]
    \end{gathered} \\
    & \text{subject to}
    & & \sum_{\mathclap{(j, i) \in \delta^{k -}_{i}}} (q_{ji}^{k+} - q_{ji}^{k-}) - \sum_{\mathclap{(i, j) \in \delta^{k +}_{i}}} (q_{ij}^{k+} - q_{ij}^{k-}) = -\overline{q}_{i}^{k}, \, \forall i \in \mathcal{D} \\
    & & & q_{ij}^{k \pm} \geq 0, \, \forall (i, j) \in \mathcal{A} \cup \mathcal{V}^{k}_{1} \cup \mathcal{P}^{k}_{1}, \, q_{ij}^{k-} = 0, \, \forall (i, j) \in \mathcal{P}^{k}_{1}.
\end{aligned}
\label{equation:p}%
\end{equation}
where $\mathcal{P}_{1}^{k} \subset \mathcal{P}$ and $\mathcal{V}_{1}^{k} \subset \mathcal{V}$ denote the sets of pumps and valves that are \emph{active} at time $k \in \mathcal{K}$, respectively (i.e., all inactive components are assumed to be ``removed'').
Similarly, $\delta_{i}^{k +}$ and $\delta_{i}^{k - }$ include only the \emph{active} set of components at time $k \in \mathcal{K}$.
Note that the second sum of the objective function models the nonlinear head gain behavior associated with active pumps.
Similar treatments of head gain appear in the experiments of \citet{collins1978solving}, the extension of the Content Model to pump networks in \citet{collins1979multiple}, and the hydraulic analysis study of \citet{moosavian2015particle}.
We emphasize that all studies stress the importance of head gain functions that are \emph{strictly concave} to ensure solution uniqueness for \eqref{equation:p}.

\section{The Co-Content Model}
\label{section:appendix-dual-derivation}
We next derive a ``Co-Content'' variant of the Content Problem \eqref{equation:p} (i.e., its dual).
This is ultimately used to derive the valid inequalities in Section \ref{section:feasibility-step}.
First, note that the objective function of Problem \eqref{equation:p} is convex and all constraints are affine.
The linearity constraint qualification thus guarantees the existence of a strong dual.
One method that can be used to derive the dual problem of $(\textnormal{P}^{k})$ is Lagrangian duality, i.e.,
\begin{equation}
    \max_{h} \min_{q^{\pm} \in \mathcal{X}} \mathcal{L}(q^{+}, q^{-}, h) = \max_{h} \ell(h) \label{eqn:lagrangian-duality},
\end{equation}
where $\mathcal{X}$ includes nonnegativity constraints, $\mathcal{L}$ is the Lagrangian of $(\textnormal{P}^{k})$ with dual variables $h$ (for flow conservation), and $\ell(h)$ is the Lagrangian dual function to be maximized.
Following $(\textnormal{P}^{k})$, its Lagrangian is
\begin{equation}
\begin{gathered}
    \mathcal{L}(q^{+}, q^{-}, h) := \sum_{i \in \mathcal{D}} h_{i}^{k} \overline{q}_{i}^{k}
    + \sum_{\mathclap{(i, j) \in \mathcal{V}^{k}_{1}}} (h_{i}^{k} - h_{j}^{k}) (q_{ij}^{k -} - q_{ij}^{k +})
    + \sum_{(i, j) \in \mathcal{A}} \left[\frac{L_{ij} r_{ij}}{1 + \alpha} (q_{ij}^{k +})^{1 + \alpha} - (h_{i}^{k} - h_{j}^{k}) q_{ij}^{k +} \right] \\
    + \sum_{(i, j) \in \mathcal{A}} \left[\frac{L_{ij} r_{ij}}{1 + \alpha} (q_{ij}^{k -})^{1 + \alpha} + (h_{i}^{k} - h_{j}^{k}) q_{ij}^{k -} \right]
    - \sum_{(i, j) \in \mathcal{P}^{k}_{1}} \left[a_{ij} q_{ij}^{k +} + \frac{b_{ij} (q_{ij}^{k +})^{c_{ij} + 1}}{c_{ij} + 1} + (h_{i}^{k} - h_{j}^{k}) q_{ij}^{k +} \right].
\end{gathered}
\label{eqn:lagrangian}
\end{equation}
First, observe that if any $(h_{i}^{k} - h_{j}^{k}) \neq 0$ for $(i, j) \in \mathcal{V}_{1}^{k}$, minimization of $\mathcal{L}(q^{+}, q^{-}, h)$ becomes unbounded, as either $q_{ij}^{k+}$ or $q_{ij}^{k-}$ could be made arbitrarily large.
Otherwise, when $h_{i}^{k} - h_{j}^{k} = 0$ for all $(i, j) \in \mathcal{V}_{1}^{k}$, Equation \eqref{eqn:lagrangian} is highly separable in the remaining $q_{ij}^{k \pm}$.
As such, minimization over $q^{\pm}$ in Equation \eqref{eqn:lagrangian-duality} is straightforward.

To derive the minimum, it suffices to minimize each component of the third through fifth sums over their corresponding $q_{ij}^{k \pm}$ while imposing nonnegativity on $q_{ij}^{k \pm}$.
Note that terms of the third and fourth sums are of the form $\frac{L_{ij} r_{ij}}{1 + \alpha} \left(q_{ij}^{k \pm}\right)^{1 + \alpha} + t_{ij}^{k} q_{ij}^{k \pm}$, where the sign of $t_{ij}^{k}$ is unknown.
There are two possibilities: if $t_{ij}^{k} \geq 0$, the component is nondecreasing in $q_{ij}^{k \pm}$ over $q_{ij}^{k \pm} \geq 0$, which implies its minimum is attained at $q_{ij}^{k \pm} = 0$.
Otherwise, if $t_{ij}^{k} < 0$, the function is decreasing at $q_{ij}^{k \pm} = 0$, attains its minimum, then starts increasing.
This minimum is attained at $\hat{q}_{ij}^{k \pm} = \sqrt[\alpha]{-\frac{t_{ij}^{k}}{L_{ij} r_{ij}}}$.
After reduction, the minimum value of the corresponding component is thus
\begin{equation}
    \left(\frac{L_{ij} r_{ij}}{1 + \alpha} - L_{ij} r_{ij}\right) \left(-\frac{t_{ij}^{k}}{L_{ij} r_{ij}}\right)^{1 + \frac{1}{\alpha}} = \frac{-\alpha}{1 + \alpha} \frac{\left(-t_{ij}^{k}\right)^{1 + \frac{1}{\alpha}}}{\sqrt[\alpha]{L_{ij} r_{ij}}} \label{eqn:dual-min}.
\end{equation}
Next, note that the third and fourth sums in Equation \eqref{eqn:lagrangian} can be paired such that the $L_{ij} r_{ij}$ coefficients of each term are the same, while the $t_{ij}^{k}$ coefficients are opposite.
That is, one term (with nonnegative $t_{ij}^{k}$) has a minimum at zero, while the other has a minimum equal to the right-hand side of Equation \eqref{eqn:dual-min}.

We next maximize terms of the final summation involving active pumps, since the summation is negated.
Note that here, all terms are of the form $f(q_{ij}^{k +}) + t_{ij}^{k} q_{ij}^{k +}$, where $f(q_{ij}^{k +})$ is a strictly concave function.
There are two possibilities: if $t_{ij}^{k} \geq 0$, the function is increasing at $q_{ij}^{k +} = 0$, attains its maximum, then starts decreasing.
Otherwise, if $t_{ij}^{k} < 0$, the function is nonincreasing in $q_{ij}^{k +}$ over $q_{ij}^{k +} \geq 0$ and thus attains its maximum at $q_{ij}^{k +} = 0$.
Analyzing stationarity conditions of the first case, $t_{ij}^{k} \geq 0$, the maximum is attained at the point
\begin{equation}
    \hat{q}_{ij}^{k +} = \left(-\frac{t_{ij}^{k} + a_{ij}}{b_{ij}}\right)^{\frac{1}{c_{ij}}}.
\end{equation}
After simplifying, then, the maximum value of the corresponding component (i.e., a term in the last sum) is
\begin{equation}
    \frac{b_{ij} \left(-\frac{t_{ij}^{k} + a_{ij}}{b_{ij}}\right)^{1 + \frac{1}{c_{ij}}}}{c_{ij}+1} + (a_{ij} + t_{ij}^{k})
    \left(-\frac{t_{ij}^{k} + a_{ij}}{b_{ij}}\right)^{\frac{1}{c_{ij}}}
\end{equation}
With the knowledge that $t_{ij}^{k} \geq 0$, $a_{ij} > 0$, $b_{ij} < 0$, and $c_{ij} > 1$, both terms are concave in $t_{ij}^{k}$, as expected.

We have determined the conditions at which all relevant terms involving $q_{ij}^{k \pm}$ are minimized in Equation \eqref{eqn:lagrangian}.
Specifically, there are two possibilities.
The first is when any $h_{i} - h_{j} \neq 0$ for $(i, j) \in \mathcal{V}_{1}^{k}$.
As discussed, this is an unbounded problem that yields $-\infty$.
The second is when $h_{i} - h_{j} = 0$ for all $(i, j) \in \mathcal{V}_{1}^{k}$, which yields
\begin{equation}
\begin{gathered}
    \tilde{\ell}(h) =
    \sum_{\mathclap{i \in \mathcal{D}}} h_{i}^{k} \overline{q}_{i}^{k} -
    \sum_{(i, j) \in \mathcal{A}} \left[\frac{\alpha}{1 + \alpha} \frac{\lvert h_{i} - h_{j} \rvert^{1 + \frac{1}{\alpha}}}{\sqrt[\alpha]{L_{ij} r_{ij}}}\right] \\
    - \sum_{(i, j) \in \mathcal{P}_{1}^{k}}
    \left[
    \frac{b_{ij} \left(\frac{h_{j}^{k} - h_{i}^{k} - a_{ij}}{b_{ij}}\right)^{1 + \frac{1}{c_{ij}}}}{c_{ij}+1}+(a_{ij} - h_{j}^{k} + h_{i}^{k})
    \left(\frac{h_{j}^{k} - h_{i}^{k} - a_{ij}}{b_{ij}}\right)^{\frac{1}{c_{ij}}}
    \right].
\end{gathered}
\end{equation}

The Lagrangian dual function we are interested in maximizing to obtain the dual \emph{problem} is thus
\begin{equation}
    \ell(h) = \begin{cases}
        \tilde{\ell}(h) & h_{i}^{k} - h_{j}^{k} = 0, \, \forall (i, j) \in \mathcal{V}_{1}^{k} \\
        -\infty & \textrm{otherwise}
    \end{cases}%
\end{equation}
Clearly, $\ell(h)$ will be maximized when the first case is satisfied.
This yields the dual maximization problem
\begin{equation}
\begin{aligned}
    & \text{maximize}
    & & \tilde{\ell}(h) \\
    & \text{subject to}
    & &  h_{i}^{k} - h_{j}^{k} = 0, \, \forall (i, j) \in \mathcal{V}_{1}^{k}.
\end{aligned}
\label{equation:d-1}%
\end{equation}
Note that since $\tilde{\ell}(h)$ is concave in head differences, Problem \eqref{equation:d-1} is also concave, as expected.
Using a standard rewriting of absolute value terms and noting that $h_{j}^{k} - h_{i}^{k} = g_{ij}^{k}$ for $(i, j) \in \mathcal{P}_{1}^{k}$ then yields
\begin{equation}
\tag{D$^{k}$}
\begin{aligned}
    & \text{maximize}
    & & \begin{gathered}
      \sum_{i \in \mathcal{D}} h_{i} \overline{q}_{i}^{k} - \frac{\alpha}{1 + \alpha} \sum_{(i, j) \in \mathcal{A}} \frac{1}{\sqrt[\alpha]{L_{ij} r_{ij}}} \left[(\Delta h^{k +}_{ij})^{1 + \frac{1}{\alpha}} + (\Delta h^{k -}_{ij})^{1 + \frac{1}{\alpha}}\right] \\
      - \sum_{(i, j) \in \mathcal{P}_{1}^{k}}
        \left[
        \frac{b_{ij} \left(\frac{g_{ij}^{k} - a_{ij}}{b_{ij}}\right)^{1 + \frac{1}{c_{ij}}}}{c_{ij}+1}+(a_{ij} - g_{ij}^{k})
        \left(\frac{g_{ij}^{k} - a_{ij}}{b_{ij}}\right)^{\frac{1}{c_{ij}}}
        \right] \\
    \end{gathered} \\
    & \text{subject to}
    & & h_{i}^{k} - h_{j}^{k} = 0, \, \forall (i, j) \in \mathcal{V}_{1}^{k}, ~ h_{j}^{k} - h_{i}^{k} = g_{ij}^{k}, \, \forall (i, j) \in \mathcal{P}_{1}^{k} \\
    & & & \Delta h_{ij}^{k +} - \Delta h_{ij}^{k -} = h_{i}^{k} - h_{j}^{k}, \, \forall (i, j) \in \mathcal{A} \\
    & & & \Delta h_{ij}^{k +}, \Delta h_{ij}^{k -} \geq 0, \, \forall (i, j) \in \mathcal{A}, \; g_{ij}^{k} \geq 0, \, \forall (i, j) \in \mathcal{P}_{1}^{k}.
\end{aligned}
\label{equation:d}%
\end{equation}
This is the dual problem that is used to derive the valid inequalities in the subsequent appendix.

\section{Duality-based Valid Convex Inequalities}
\label{section:micp_reformulation}
Using the primal-dual pair of problems, \eqref{equation:p} and \eqref{equation:d}, in this appendix, we derive novel duality-based valid inequalities for operational WDN problems, including the OWF.
Similar to \citet{tasseff2020exact}, letting the objective of Problem \eqref{equation:p} be denoted by $f_{P}^{k}(q^{k})$ and the objective of Problem \eqref{equation:d} be denoted by $f_{D}^{k}(h^{k})$, the \emph{convex} strong duality constraint $f_{P}^{k}(q^{k}) - f_{D}^{k}(h^{k}) \leq 0$ reformulates the nonconvex physics at time index $k \in \mathcal{K}$.
Noting that $\Delta h_{ij}^{k \pm} = L_{ij} r_{ij} (q_{ij}^{k \pm})^{\alpha}$ for pipes and $g_{ij}^{k} = a_{ij} + b_{ij} (q_{ij}^{k +})^{c_{ij}}$ for pumps, this simplifies to
\begin{equation}
    \begin{gathered}
        \sum_{\mathclap{(i, j) \in \mathcal{A}}} L_{ij} r_{ij} \left[(q_{ij}^{k +})^{1 + \alpha} + (q_{ij}^{k -})^{1 + \alpha}\right]
        - \sum_{\mathclap{(i, j) \in \mathcal{P}}} \left[a_{ij} q_{ij}^{k +} + b_{ij} (q_{ij}^{k +})^{c_{ij} + 1}\right] \\
        - \sum_{i \in \mathcal{R} \cup \mathcal{T}} h_{i}^{k} \left[\sum_{(i, j) \in \delta^{k +}_{i}} (q_{ij}^{k +} - q_{ij}^{k -}) - \sum_{(j, i) \in \delta^{k -}_{i}} (q_{ji}^{k +} - q_{ji}^{k -})\right]
       - \sum_{i \in \mathcal{D}} h_{i} \overline{q}_{i}^{k} \leq 0.
    \end{gathered}
\end{equation}
There is an important caveat, however: the convexity of this inequality assumes that tank heads are fixed, along with the conventional assumption that all reservoir heads and demand flows are also fixed.
That is, in the context of the OWF, unlike the optimal WDN design problem of \citet{tasseff2020exact}, the strong duality constraint cannot be directly embedded to provide an \emph{exact} convex reformulation of the system physics.

Unlike the inequalities of \citet{tasseff2020exact}, the operational constraints we derive include products of head and flow variables at tanks.
These bilinear terms must be relaxed to ensure convexity.
This is accomplished by using a standard McCormick relaxation of the complicating terms.
That is, products of flows and heads, denoted here as $q_{i}^{k} h_{i}^{k}$ for brevity, are relaxed in the inequalities by introducing the constraints
\begin{subequations}
\begin{align}
    w_{i}^{k} &\geq \underline{q}_{i}^{k} h_{i}^{k} + \underline{h}_{i}^{k} q_{i}^{k} - \underline{q}_{i}^{k} \underline{h}_{i}^{k}, ~ \forall i \in \mathcal{T}, ~ \forall k \in \mathcal{K} \\
    w_{i}^{k} &\geq \overline{q}_{i}^{k} h_{i}^{k} + \overline{h}_{i}^{k} q_{i}^{k} - \overline{q}_{i}^{k} \overline{h}_{i}^{k}, ~ \forall i \in \mathcal{T}, ~ \forall k \in \mathcal{K} \\
    w_{i}^{k} &\leq \underline{q}_{i}^{k} h_{i}^{k} + \overline{h}_{i}^{k} q_{i}^{k} - \underline{q}_{i}^{k} \overline{h}_{i}^{k}, ~ \forall i \in \mathcal{T}, ~ \forall k \in \mathcal{K} \\
    w_{i}^{k} &\leq \overline{q}_{i}^{k} h_{i}^{k} + \underline{h}_{i}^{k} q_{i}^{k} - \overline{q}_{i}^{k} \underline{h}_{i}^{k}, ~ \forall i \in \mathcal{T}, ~ \forall k \in \mathcal{K},
\end{align}
\label{equation:mccormick-convex}%
\end{subequations}
where $w_{i}^{k}$, $i \in \mathcal{T}$, $k \in \mathcal{K}$, denotes the product of tank head and flow.
The inequalities are then rewritten as
\begin{equation}
    \begin{gathered}
        \sum_{\mathclap{(i, j) \in \mathcal{A}}} L_{ij} r_{ij} \left[(q_{ij}^{k +})^{1 + \alpha} + (q_{ij}^{k -})^{1 + \alpha}\right]
       - \sum_{\mathclap{(i, j) \in \mathcal{P}}} \left[a_{ij} q_{ij}^{k +} + b_{ij} (q_{ij}^{k +})^{c_{ij} + 1}\right]
       - \sum_{\mathclap{i \in \mathcal{T}}} w_{i}^{k} \\
       - \sum_{i \in \mathcal{R}} h_{i}^{k} \left[\sum_{(i, j) \in \delta^{k +}_{i}} (q_{ij}^{k +} - q_{ij}^{k -}) - \sum_{(j, i) \in \delta^{k -}_{i}} (q_{ji}^{k +} - q_{ji}^{k -})\right]
       - \sum_{i \in \mathcal{D}} h_{i} \overline{q}_{i}^{k} \leq 0, \, \forall k \in \mathcal{K}.
    \end{gathered}
\label{equation:convex_inequality-2}%
\end{equation}
Note that \emph{these} constraints are convex and can be applied for all time indices $k \in \mathcal{K}$.
Also, observe that the strength of each inequality depends on the tightness of tank head and flow bounds, as used in the McCormick relaxation Constraints \eqref{equation:mccormick-convex}.
We remark, however, that the range of tank heads is often small and limited by predefined tank level bounds.
This concludes the derivation of the inequalities presented in Section \ref{section:feasibility-step}.

\section{Additional Optimization-based Valid Inequalities}
\begin{algorithm}[t]
\caption{Binary-binary OBCG procedure, where one variable is fixed to one.}
\label{algorithm:obcg-2}
\begin{algorithmic}[1]
  \Input Any direction-based model formulation, which comprises some feasible set $\Omega$
  \Output Valid inequalities for direction-based model formulations, denoted by $\bar{\mathcal{X}}$
  \State
  $\mathcal{B} \gets \{y_{ij}^{k} : (i, j) \in \mathcal{L}\} \cup \{z_{ij}^{k} : (i, j) \in \mathcal{P} \cup \mathcal{V}\}$, $\bar{\mathcal{X}} = \emptyset$ \label{line:obcg-2-collect-variable-sets}
  \ForAll{$(x_{1}, x_{2}) \in (\mathcal{B} \times \mathcal{B}) \setminus \{(x, x) : x \in \mathcal{B}\}$} \label{line:obcg-2-binary-binary-loop}
      \State $(\underline{x}_{1}^{1}, \overline{x}_{1}^{1}) \gets (\textrm{minimize} ~ x_{1} ~ \textrm{s.t.} ~ \Omega \cup \{x_{2} = 1\}, ~ \textrm{maximize} ~ x_{1} ~ \textrm{s.t.} ~ \Omega \cup \{x_{2} = 1\})$ \label{line:obcg-2-binary-binary-minimize-1}
      \If{$\underline{x}_{1}^{1} = \overline{x}_{1}^{1} = 0$} \label{line:obcg-2-if-binary-binary-0-0}
          \State $\bar{\mathcal{X}} \gets \bar{\mathcal{X}} \cup \{x_{1} + x_{2} \leq 1\}$ \label{line:obcg-2-cut-binary-binary-0-0}
      \ElsIf{$\underline{x}_{1}^{1} = \overline{x}_{1}^{1} = 1$} \label{line:obcg-2-if-binary-binary-0-1}
          \State $\bar{\mathcal{X}} \gets \bar{\mathcal{X}} \cup \{x_{1} \geq x_{2}\}$ \label{line:obcg-2-cut-binary-binary-0-1}
      \EndIf \label{line:obcg-2-endif-binary-binary-0}
  \EndFor
\end{algorithmic}
\end{algorithm}

The OBCG Algorithm \ref{algorithm:obcg-2} is an analogue of the OBCG Algorithm \ref{algorithm:obcg-1}, where the second variable is fixed to one instead of zero.
As in Algorithm \ref{algorithm:obcg-1}, the goal of Algorithm \ref{algorithm:obcg-2} is to derive cuts that improve continuous variable bounds \emph{depending on the value of each discrete variable}.
Here, Line \ref{line:obcg-2-collect-variable-sets} instantiates the set of binary variables to consider in the cut generation procedure.
Line \ref{line:obcg-2-binary-binary-loop} defines the binary-binary cut generation loop, where all unique, ordered pairs of binary variables are considered.
Line \ref{line:obcg-2-binary-binary-minimize-1} minimizes and maximizes the first variable in the pair, respectively, subject to the relaxation constraints and a fixing of the second variable to one.
On Line \ref{line:obcg-2-if-binary-binary-0-0}, if both minimization and maximization imply objective values of zero, a cut relating the two discrete variables can be derived, which is added to the set of cuts $\bar{\mathcal{X}}$ on Line \ref{line:obcg-2-cut-binary-binary-0-0}.
On the other hand, if minimization and maximization both yield objective values of one, a similar process is followed on Lines \ref{line:obcg-2-if-binary-binary-0-1}--\ref{line:obcg-2-endif-binary-binary-0}.

\section{Network Graph Visualizations}
In this subsection, graph visualizations of the three test networks are provided.
Simple FSD is illustrated in Figure \ref{fig:graph_simple_fsd}, AT(M) in Figure \ref{fig:graph_atm}, and Poormond in Figure \ref{fig:graph_poormond}.
Node and link labels correspond, first, to the original data set label (or component name), followed by the \textsc{WaterModels}-specific integer index in parentheses.
The elevations of the nodes are color-coded according the colorbar shown above each network visualization, where the elevation is in meters.
Demand nodes are indicated by ovals, and nodes with nonzero demands have the volumetric flow rate (in m$^{3}$/s) of the demand at time index $k = 1$ included textually in the label.
Reservoirs are indicated by diamonds and tanks by squares.
Pipes are shown as black lines, where the thickness indicates the relative pipe diameter (i.e., thicker lines correspond to larger diameter pipes).
Pipe lengths are noted in the pipe labels.
Pumps are shown as red lines, and check valves as blue lines.
Arrows indicate links that are unidirectional, while links without arrows may have flow in either direction.

\begin{figure}[t]
\centering
\includegraphics[width=0.5\textwidth]{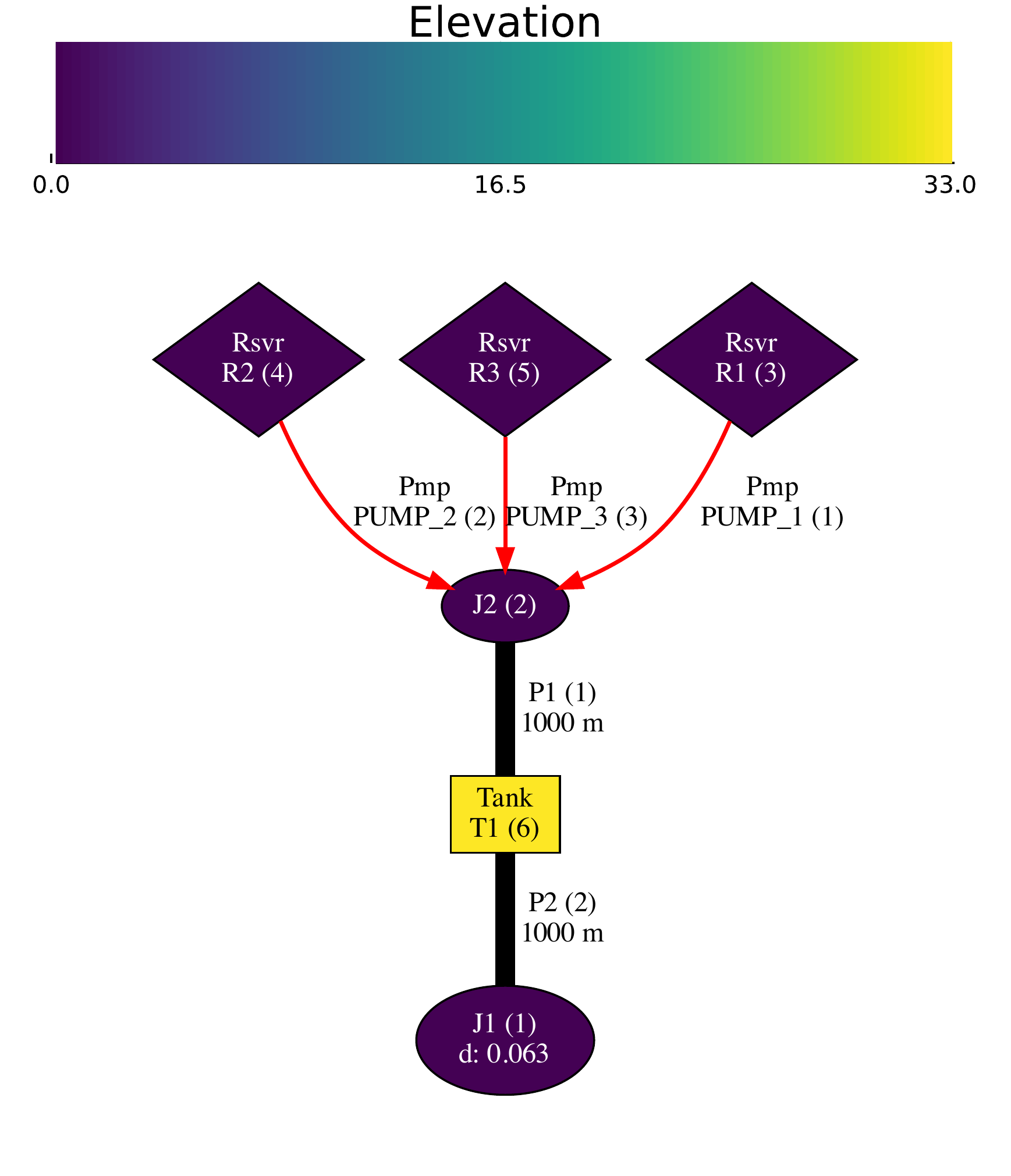}
\caption{Visualization of the Simple FSD water network.}
\label{fig:graph_simple_fsd}
\end{figure}

\clearpage

\begin{figure}[t]
\centering
\includegraphics[height=\textheight]{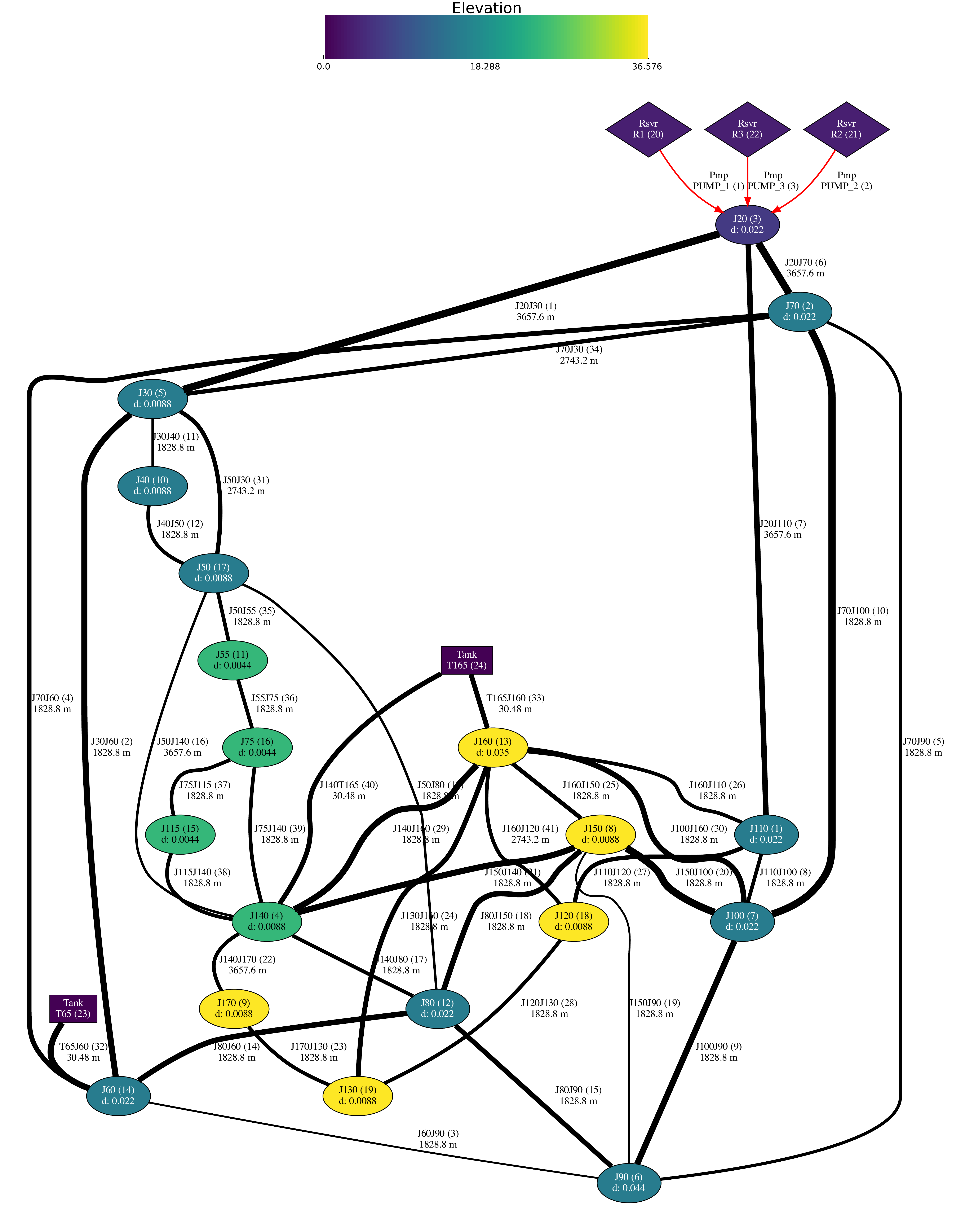}
\caption{Visualization of the AT(M) water network.}
\label{fig:graph_atm}
\end{figure}

\clearpage

\begin{figure}[t]
\centering
\includegraphics[height=\textheight]{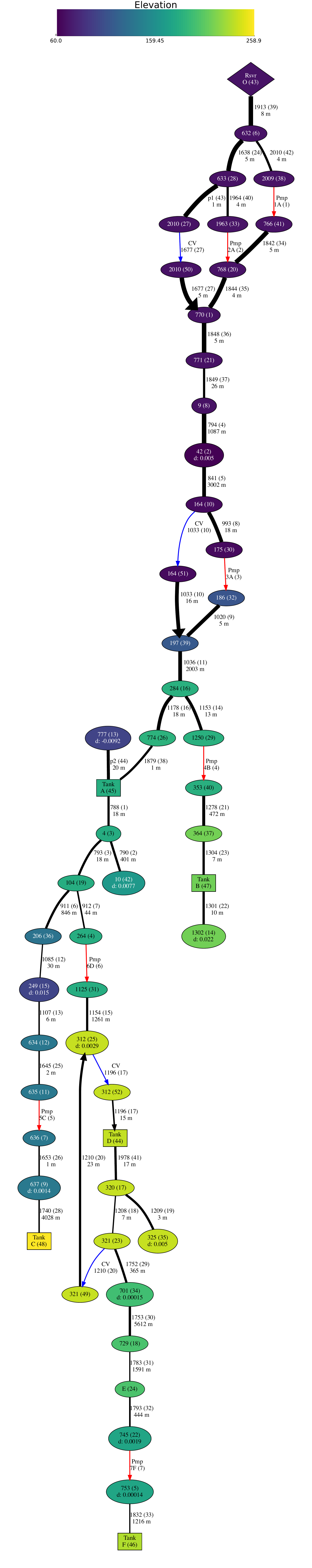}
\caption{Visualization of the Poormond water network.}
\label{fig:graph_poormond}%
\end{figure}

\clearpage

\section{Electricity Price Profiles}
Figures \ref{figure:prices-simple_fsd}, \ref{figure:prices-matm}, and \ref{figure:prices-poormond} illustrate the electricity price profiles for a subset of the OWF instances introduced in Section \ref{subsection:experimental_setup}.
Specifically, the five unique price profiles that define the five separate instances corresponding to $K = 48$ are depicted for each of the three networks.
Price profiles for $K = 12$ and $K = 24$ can be thought of as coarser representations of these $K = 48$ price profiles and are omitted from these plots for clarity.

\begin{figure}[t]
    \centering
    \begin{tikzpicture}
	\centering
	\begin{axis}[legend cell align=left,enlargelimits=false,xtick pos=left,ytick pos=left,
                 legend columns=1,xlabel=Time since start (seconds),ylabel=Electricity cost per Joule (\textsterling),
                 height=6.0cm,width=\textwidth,/pgf/number format/.cd,1000 sep={},
                 legend style={at={(0.02, 0.98)},anchor=north west}]


        \pgfplotstableread[col sep = comma]{Simple_FSD-48_Steps-Day_21-Electricity_Prices.csv}{\data};
        \addplot[very thick, mark=*, opacity=1.0, index of colormap=0 of Dark2-8] table [x = time, y = energy_price]{\data};
        \addlegendentry{$K = 48$, Day 1}



        \pgfplotstableread[col sep = comma]{Simple_FSD-48_Steps-Day_22-Electricity_Prices.csv}{\data};
        \addplot[very thick, mark=*, opacity=1.0, index of colormap=1 of Dark2-8] table [x = time, y = energy_price]{\data};
        \addlegendentry{$K = 48$, Day 2}



        \pgfplotstableread[col sep = comma]{Simple_FSD-48_Steps-Day_23-Electricity_Prices.csv}{\data};
        \addplot[very thick, mark=*, opacity=1.0, index of colormap=2 of Dark2-8] table [x = time, y = energy_price]{\data};
        \addlegendentry{$K = 48$, Day 3}



        \pgfplotstableread[col sep = comma]{Simple_FSD-48_Steps-Day_24-Electricity_Prices.csv}{\data};
        \addplot[very thick, mark=*, opacity=1.0, index of colormap=3 of Dark2-8] table [x = time, y = energy_price]{\data};
        \addlegendentry{$K = 48$, Day 4}



        \pgfplotstableread[col sep = comma]{Simple_FSD-48_Steps-Day_25-Electricity_Prices.csv}{\data};
        \addplot[very thick, mark=*, opacity=1.0, index of colormap=4 of Dark2-8] table [x = time, y = energy_price]{\data};
        \addlegendentry{$K = 48$, Day 5}




	\end{axis}
\end{tikzpicture}
    \caption{Five price profiles for the five $K = 48$ Simple FSD OWF instances considered in the paper.}
    \label{figure:prices-simple_fsd}%
\end{figure}
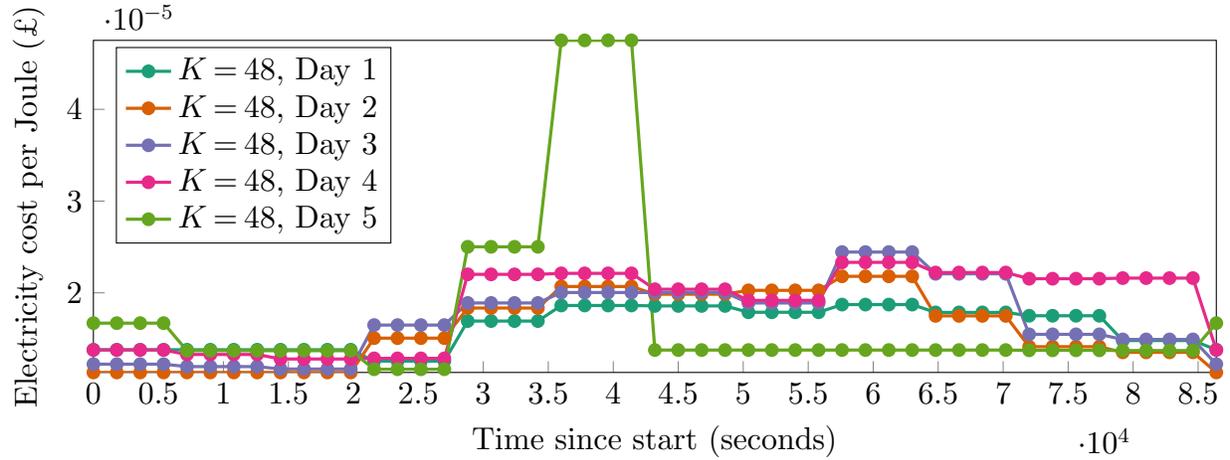

\begin{figure}[t]
    \centering
    \begin{tikzpicture}
	\centering
	\begin{axis}[legend cell align=left,enlargelimits=false,xtick pos=left,ytick pos=left,
                 legend columns=1,xlabel=Time since start (seconds),ylabel=Electricity cost per Joule (\textsterling),
                 height=6.0cm,width=\textwidth,/pgf/number format/.cd,1000 sep={},
                 legend style={at={(0.02, 0.98)},anchor=north west}]


        \pgfplotstableread[col sep = comma]{MATM-48_Steps-Day_21-Electricity_Prices.csv}{\data};
        \addplot[very thick, mark=*, opacity=1.0, index of colormap=0 of Dark2-8] table [x = time, y = energy_price]{\data};
        \addlegendentry{$K = 48$, Day 1}



        \pgfplotstableread[col sep = comma]{MATM-48_Steps-Day_22-Electricity_Prices.csv}{\data};
        \addplot[very thick, mark=*, opacity=1.0, index of colormap=1 of Dark2-8] table [x = time, y = energy_price]{\data};
        \addlegendentry{$K = 48$, Day 2}



        \pgfplotstableread[col sep = comma]{MATM-48_Steps-Day_23-Electricity_Prices.csv}{\data};
        \addplot[very thick, mark=*, opacity=1.0, index of colormap=2 of Dark2-8] table [x = time, y = energy_price]{\data};
        \addlegendentry{$K = 48$, Day 3}



        \pgfplotstableread[col sep = comma]{MATM-48_Steps-Day_24-Electricity_Prices.csv}{\data};
        \addplot[very thick, mark=*, opacity=1.0, index of colormap=3 of Dark2-8] table [x = time, y = energy_price]{\data};
        \addlegendentry{$K = 48$, Day 4}



        \pgfplotstableread[col sep = comma]{MATM-48_Steps-Day_25-Electricity_Prices.csv}{\data};
        \addplot[very thick, mark=*, opacity=1.0, index of colormap=4 of Dark2-8] table [x = time, y = energy_price]{\data};
        \addlegendentry{$K = 48$, Day 5}




	\end{axis}
\end{tikzpicture}
    \caption{Five price profiles for the five $K = 48$ AT(M) OWF instances considered in the paper.}
    \label{figure:prices-matm}%
\end{figure}
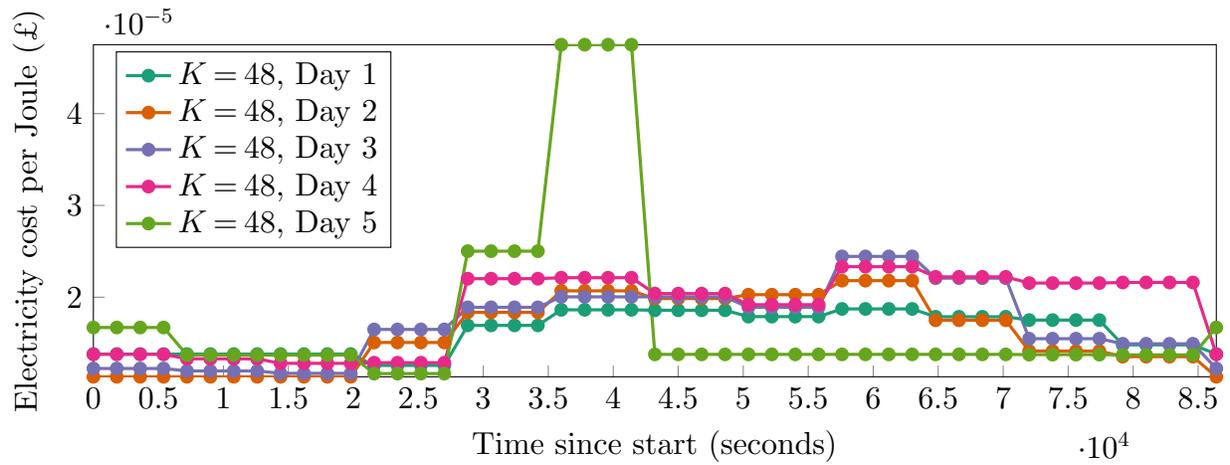

\begin{figure}[t]
    \centering
    \begin{tikzpicture}
	\centering
	\begin{axis}[legend cell align=left,enlargelimits=false,xtick pos=left,ytick pos=left,
                 legend columns=1,xlabel=Time since start (seconds),ylabel=Electricity cost per Joule (\textsterling),
                 height=6.0cm,width=\textwidth,/pgf/number format/.cd,1000 sep={},
                 legend style={at={(0.02, 0.98)},anchor=north west}]


        \pgfplotstableread[col sep = comma]{Poormond-48_Steps-Day_21-Electricity_Prices.csv}{\data};
        \addplot[very thick, mark=*, opacity=1.0, index of colormap=0 of Dark2-8] table [x = time, y = energy_price]{\data};
        \addlegendentry{$K = 48$, Day 1}



        \pgfplotstableread[col sep = comma]{Poormond-48_Steps-Day_22-Electricity_Prices.csv}{\data};
        \addplot[very thick, mark=*, opacity=1.0, index of colormap=1 of Dark2-8] table [x = time, y = energy_price]{\data};
        \addlegendentry{$K = 48$, Day 2}



        \pgfplotstableread[col sep = comma]{Poormond-48_Steps-Day_23-Electricity_Prices.csv}{\data};
        \addplot[very thick, mark=*, opacity=1.0, index of colormap=2 of Dark2-8] table [x = time, y = energy_price]{\data};
        \addlegendentry{$K = 48$, Day 3}



        \pgfplotstableread[col sep = comma]{Poormond-48_Steps-Day_24-Electricity_Prices.csv}{\data};
        \addplot[very thick, mark=*, opacity=1.0, index of colormap=3 of Dark2-8] table [x = time, y = energy_price]{\data};
        \addlegendentry{$K = 48$, Day 4}



        \pgfplotstableread[col sep = comma]{Poormond-48_Steps-Day_25-Electricity_Prices.csv}{\data};
        \addplot[very thick, mark=*, opacity=1.0, index of colormap=4 of Dark2-8] table [x = time, y = energy_price]{\data};
        \addlegendentry{$K = 48$, Day 5}




	\end{axis}
\end{tikzpicture}
    \caption{Five price profiles for the five $K = 48$ Poormond OWF instances considered in the paper.}
    \label{figure:prices-poormond}%
\end{figure}
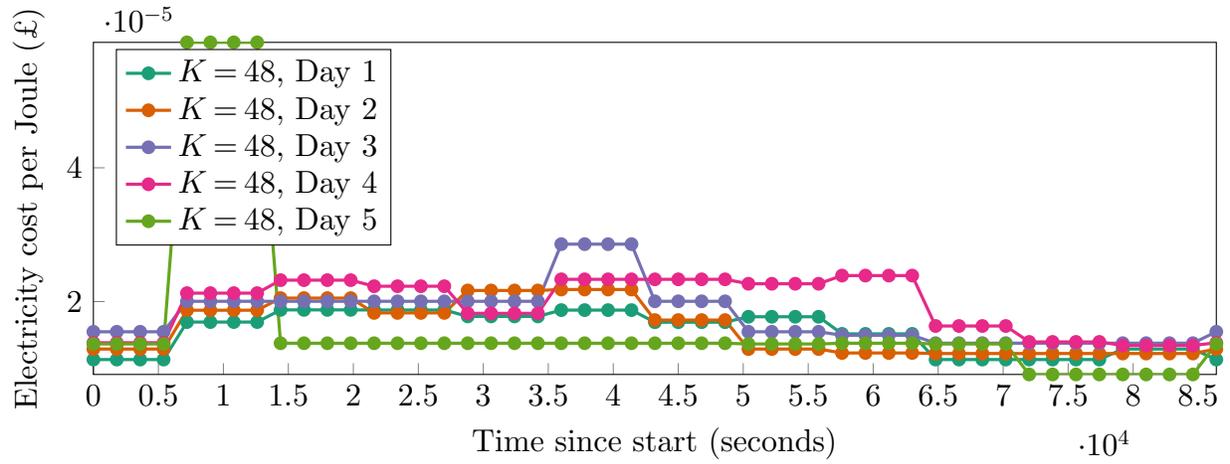

\end{APPENDICES}

\end{document}